\documentclass[11pt]{article}
\usepackage{amsmath,amsfonts,amssymb,algorithm,algorithmic,psfrag}
\usepackage{amsthm}
\usepackage{array}
\usepackage{psfrag}
\usepackage{subfigure}
\usepackage{bm}
\usepackage{color}
\usepackage[curve]{xypic}
\usepackage{bigstrut}

\newtheorem{theorem}{Theorem}[section]
\newtheorem{lemma}[theorem]{Lemma}

\theoremstyle{definition}

\theoremstyle{remark}
\newtheorem{remark}[theorem]{Remark}

\numberwithin{equation}{section}

\usepackage{xypic}

\usepackage{color} 
\usepackage{ifpdf}
\ifpdf
    \usepackage[pdftex]{graphicx}
    \usepackage[pdftex]{hyperref}
    \hypersetup{
        unicode=false,          
        pdftoolbar=true,        
        pdfmenubar=true,        
        pdffitwindow=false,     
        pdfstartview={FitH},    
        pdftitle={AKG Article},      
        pdfauthor={Andrew Gillette},   
        pdfsubject={Mathematics},    
        pdfcreator={Andrew Gillette},  
        pdfproducer={Andrew Gillette}, 
        pdfkeywords={math, mathematics}, 
        pdfnewwindow=true,      
        colorlinks=true,        
        linkcolor=red,          
        citecolor=blue,         
        filecolor=magenta,      
        urlcolor=cyan           
    }

\else
    \usepackage{graphicx}
    \usepackage{pstricks}
    
    \newcommand{\href}[2]{#2}
\fi

\usepackage[mathscr]{euscript}
\usepackage{amsbsy}

\usepackage{xypic}

\usepackage{amssymb}

\DeclareMathOperator{\Div}{div}
\DeclareMathOperator{\rot}{rot}
\DeclareMathOperator{\grad}{grad}
\DeclareMathOperator{\curl}{curl}

\DeclareMathOperator{\tr}{tr}

\renewcommand{\div}{\Div}

\newcommand{\rmH}{\mathrm{H}}

\newcommand{\bbN}{\mathbb{N}}

\newcommand{\bbR}{\mathbb{R}}

\newcommand{\calH}{\mathcal{H}}

\newcommand{\calJ}{\mathcal{J}}

\newcommand{\calP}{\mathcal{P}}
\newcommand{\calQ}{\mathcal{Q}}

\newcommand{\calS}{\mathcal{S}}

\newcommand{\bfa}{\textbf{a}}
\newcommand{\bfu}{\textbf{u}}
\newcommand{\bfv}{\textbf{v}}
\newcommand{\bfx}{\textbf{x}}

\newcommand{\dof}{\textsc{dof}}

\newcommand{\beq}{\begin{equation}}
\newcommand{\eeq}{\end{equation}}

\newcommand{\ds}{\displaystyle}

\newcommand{\bs}{{\scriptscriptstyle \bullet}}

\newcommand{\twovec}[2]{
\left[\!\!\!
\begin{array}{c}
#1\\
#2
\end{array}
\!\!\! \right]
}
\newcommand{\threevec}[3]{
\left[\!\!\!
\begin{array}{c}
#1\\
#2\\
#3
\end{array}
\!\!\! \right]
}

\newcommand{\p}{\partial}
\newcommand{\raw}{\rightarrow}

\newcommand{\ldeg}{\textnormal{ldeg}}

\newcommand{\nedelec}{N{\'e}d{\'e}lec}

\newtheorem{theo}{Th\'eor\`eme}[section]
\newtheorem{prop}[theo]{Proposition}

\theoremstyle{plain}

\title{Trimmed Serendipity Finite Element Differential Forms}

\author{Andrew Gillette\footnote{
Department of Mathematics,
University of Arizona,
Tucson, AZ 85721. \hfill
{\it agillette@math.arizona.edu}}
~and Tyler Kloefkorn\footnote{
Department of Mathematics,
University of Arizona,
Tucson, AZ 85721. \hfill
{\it tkloefkorn@math.arizona.edu}}
}

\begin{document}
\maketitle

\begin{abstract}
We introduce the family of trimmed serendipity finite element differential form spaces, defined on cubical meshes in any number of dimensions, for any polynomial degree, and for any form order.
The relation between the trimmed serendipity family and the (non-trimmed) serendipity family developed by Arnold and Awanou [\textit{Math.\ Comp.\ 83(288) 2014}] is analogous to the relation between the trimmed and (non-trimmed) polynomial finite element differential form families on simplicial meshes from finite element exterior calculus. 
We provide degrees of freedom in the general setting and prove that they are unisolvent for the trimmed serendipity spaces.
The sequence of trimmed serendipity spaces with a fixed polynomial order $r$ provides an explicit example of a system described by Christiansen and Gillette [\textit{ESAIM:M2AN 50(3) 2016}], namely, a minimal compatible finite element system on squares or cubes containing order $r-1$ polynomial differential forms.
\end{abstract}

\maketitle

\section{Introduction}
\label{sec:intro}

The `Periodic Table of the Finite Elements'~\cite{AL2014} identifies four families of polynomial differential form spaces: $\calP_r^-\Lambda^k$, $\calP_r\Lambda^k$, $\calQ_r^-\Lambda^k$ and $\calS_r\Lambda^k$.
The families $\calP_r^-\Lambda^k$ and $\calP_r\Lambda^k$ define finite element spaces on $n$-simplices while $\calQ_r^-\Lambda^k$ and $\calS_r\Lambda^k$ define finite element spaces on $n$-dimensional cubes.
In this paper, we present a fifth family, $\calS_r^-\Lambda^k$ that is closely related to but distinct from the serendipity family $\calS_r\Lambda^k$~\cite{AA2014}.
In particular, the relationships between the families $\calS_r^-\Lambda^k$ and $\calS_r\Lambda^k$ are analogous to the relationships between $\calP_r^-\Lambda^k$ and $\calP_r\Lambda^k$.

We first define the $\calS_r^-\Lambda^k$ spaces as
\begin{equation*}
\calS_r^-\Lambda^k :=\calS_{r-1}\Lambda^k +\kappa\calS_{r-1}\Lambda^{k+1},
\end{equation*}
where $\kappa$ denotes the Koszul operator.
The $\calS_r^-\Lambda^k$ spaces nest in between serendipity spaces via the inclusions:
\[\calS_{r}\Lambda^k \subset \calS_{r+1}^-\Lambda^k \subset \calS_{r+1}\Lambda^k.\]
The exterior derivative $d$ makes $\calS_r^-\Lambda^\bs$ into a cochain complex and the associated sequence
\[ 0 \rightarrow \bbR \rightarrow \calS_r^-\Lambda^0 \rightarrow \calS_r^-\Lambda^1 \rightarrow \cdots \rightarrow \calS_r^-\Lambda^{n-1} \rightarrow \calS_r^-\Lambda^n \rightarrow 0  \]
is exact.
The spaces in the above sequence have minimal dimension for $n=2$ or $n=3$ in the following sense: the sequence is a minimal compatible finite element system on $n$-cubes that contains $\calP_{r-1}\Lambda^k$ for each $k$.
All the results just mentioned, as well as others identified in this paper, hold true if $\calP$ is put in place of $\calS$ and the spaces are taken over $n$-simplices instead of $n$-cubes.
Since $\calP_r^-\Lambda^k$ spaces have been called trimmed polynomial spaces, we refer to the $\calS_r^-\Lambda^k$ spaces as trimmed serendipity spaces.

We describe the trimmed serendipity family of finite elements using the language of finite element exterior calculus~(FEEC)~\cite{AFW2006,AFW2010}.
The FEEC framework has also been used to describe the famous elements of \nedelec~\cite{N1980,N1986}, Raviart-Thomas~\cite{RT1977}, and Brezzi-Douglas-Marini~\cite{BDM85}, as well as the more recently defined elements of Arnold and Awanou~\cite{AA2011,AA2014}.
Here, we show how the FEEC framework can also describe the recently defined $AC_r$ elements on squares of Arbogast and Correa~\cite{AC2015}, the $S_{2,k}$ elements on squares and cubes of Cockburn and Fu~\cite{CF2016}, and the virtual element serendipity spaces $VEMS^f_{r,r,r-1}$ of Beir\~{a}o da Veiga, Brezzi, Marini, and Russo~\cite{BeiEtAl16}. 
A detailed comparison to these newer elements is given at the end of Section~\ref{sec:bkgd}.
Two key features of our approach that distinguish is from related papers are: ($i$) a generalized definition of degrees of freedom suitable for any $r\geq 1$, $n\geq 1$ and $0\leq k\leq n$; and ($ii$) the extensive use of tools from exterior calculus, allowing generalization to arbitrary dimension $n$ and instant coordination with other results from FEEC.

Christiansen and Gillette~\cite{CG2015} raised the question of a minimal compatible finite element system on $n$-cubes containing $\calP_{r-1}\Lambda^k$ and computed the number of degrees of freedom that such a system would need to associate to the interior of an $n$-cube, $\square_n$.
While we do not use the harmonic extension approach of~\cite{CG2015} to construct the $\calS_r^-\Lambda^k$ spaces, we do recover the expected degree of freedom counts associated to each piece of the cubical geometry.
We state the dimension of $\calS_r^-\Lambda^k(\square_n)$ for $1\leq n\leq 4$, $0\leq k\leq n$, and $1\leq r\leq 7$ in Table~\ref{tab:srm-dim}.

The spaces $\calS_r^-\Lambda^1(\square_3)$ and $\calS_r^-\Lambda^2(\square_3)$ are of potentially great interest to the computational electromagnetics community as they can be used in $H(\curl)$- and $H(\Div)$-conforming methods on meshes of affinely-mapped cubes.
Their dimensions satisfy $\dim\calS_r^-\Lambda^1(\square_3)<\dim\calS_r\Lambda^1(\square_3)$ and $\dim\calS_r^-\Lambda^2(\square_3)<\dim\calS_r\Lambda^2(\square_3)$ as well as $\dim\calS_r^-\Lambda^1(\square_3)\leq\dim\calQ_r^-\Lambda^1(\square_3)$ and $\dim\calS_r^-\Lambda^2(\square_3)\leq\dim\calQ_r^-\Lambda^2(\square_3)$, with equality only in the the case $r=1$.
Hence, a significant savings in degrees of freedom should be possible, compared to tensor product and even serendipity methods.
At the end of Section~\ref{sec:unisolv-min}, we present an example illustrating the reduction in the degrees of freedom in the context of a mixed method for Poisson's problem.

\begin{figure}[t]
\centering
\begin{equation*}
\xymatrix @R=.02in{
\parbox{.22\textwidth}
{\includegraphics[width=.22\textwidth]{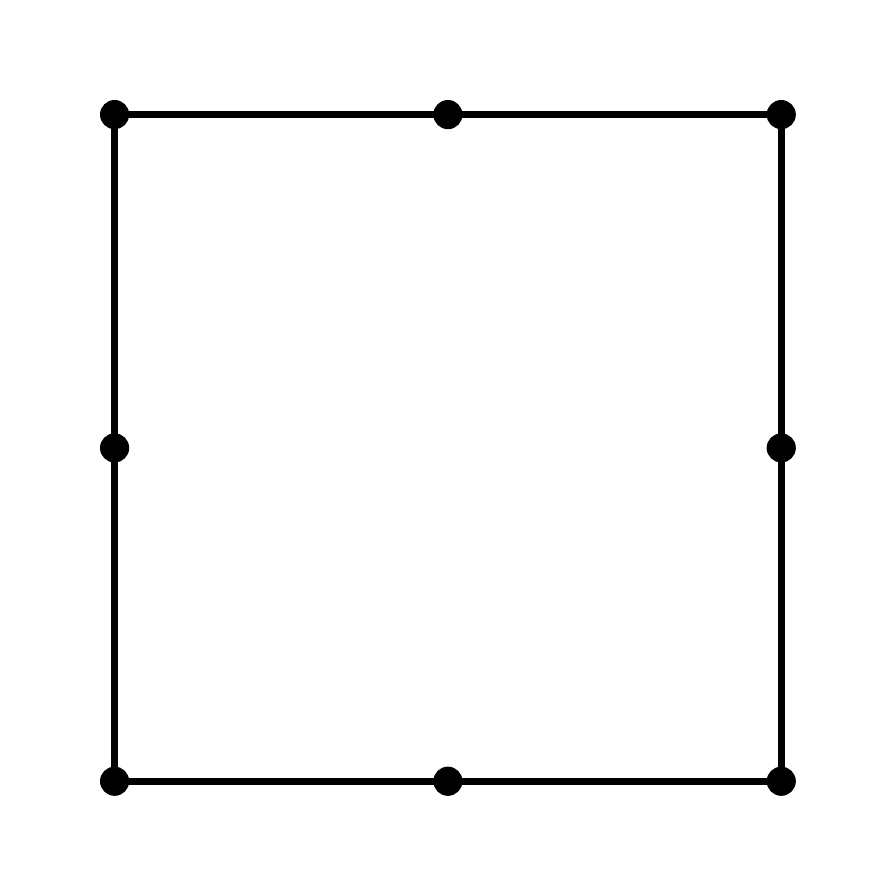}}  \ar[r]^-{d_0}_{\grad}
&
\parbox{.22\textwidth}
{\includegraphics[width=.22\textwidth]{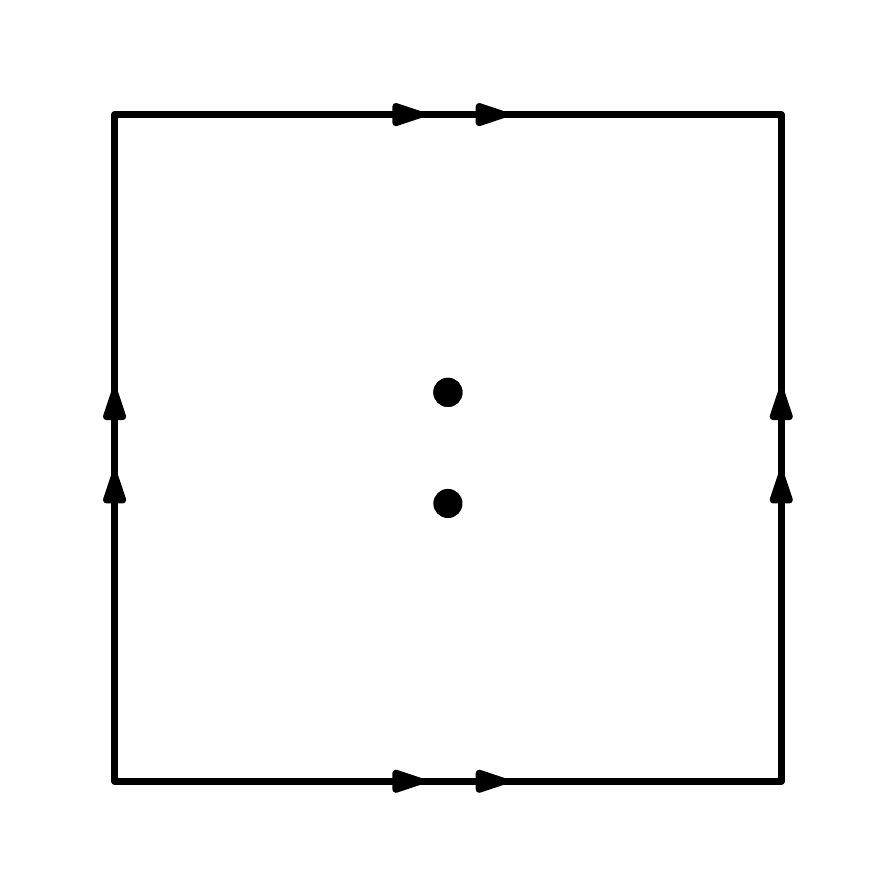}} \ar[r]^-{d_1}_{\div\rot}
&
\parbox{.22\textwidth}
{\includegraphics[width=.22\textwidth]{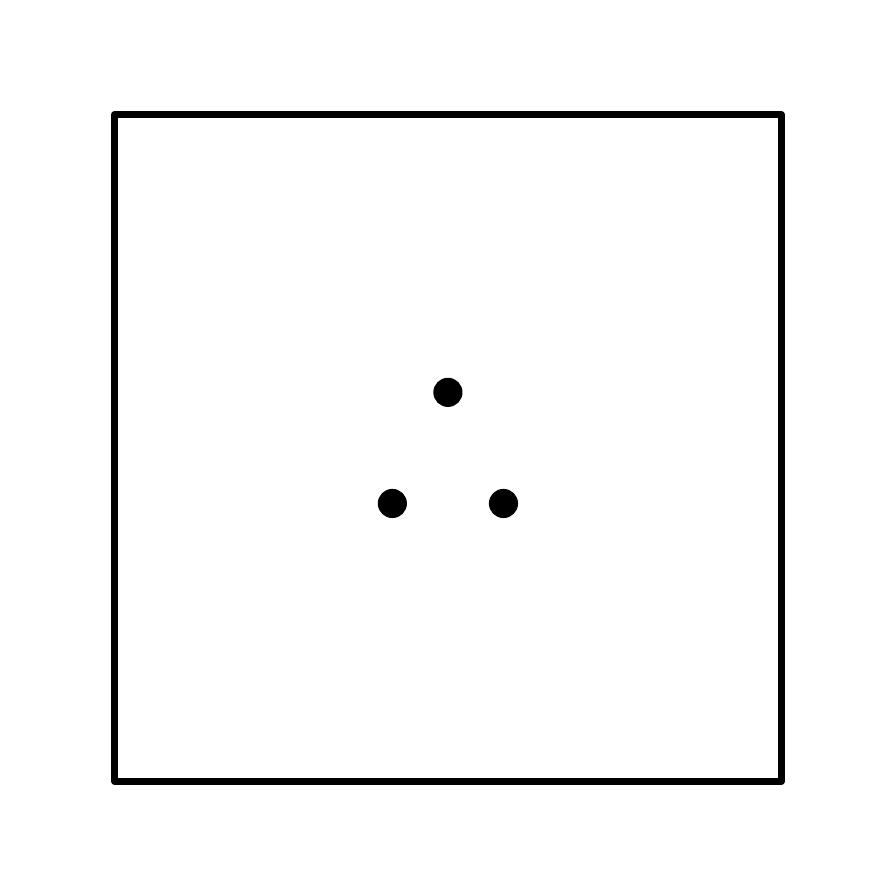}} \\
\dim \calS_2^-\Lambda^0(\square_2)=8 & \dim \calS_2^-\Lambda^1(\square_2)=10 & \dim \calS_2^-\Lambda^2(\square_2)=3 \\
}
\end{equation*}
\vspace{-2mm}
\begin{equation*}
\xymatrix @R=.02in{
\parbox{.22\textwidth}
{\includegraphics[width=.22\textwidth]{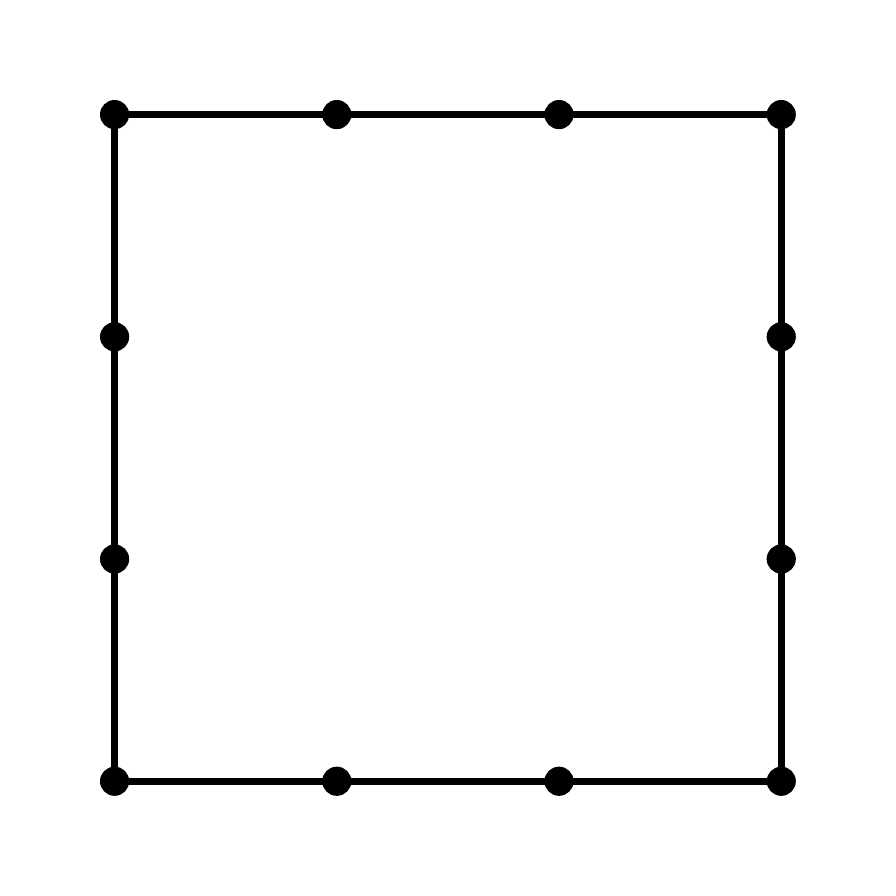}}  \ar[r]^-{d_0}_{\grad}
&
\parbox{.22\textwidth}
{\includegraphics[width=.22\textwidth]{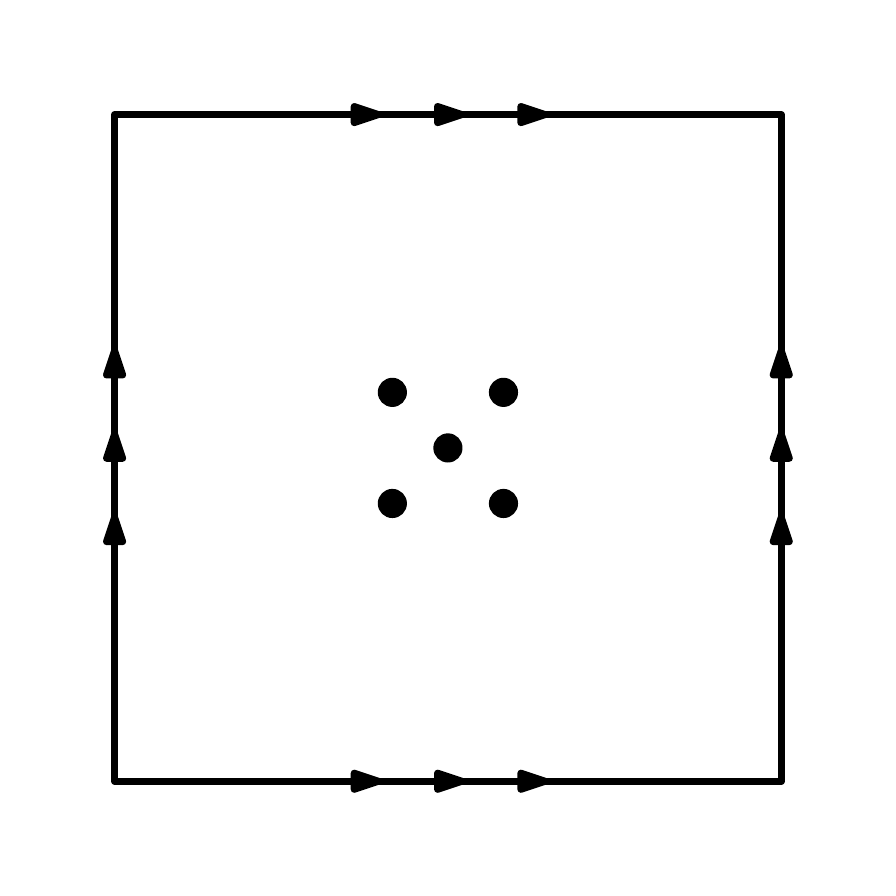}} \ar[r]^-{d_1}_{\div\rot}
&
\parbox{.22\textwidth}
{\includegraphics[width=.22\textwidth]{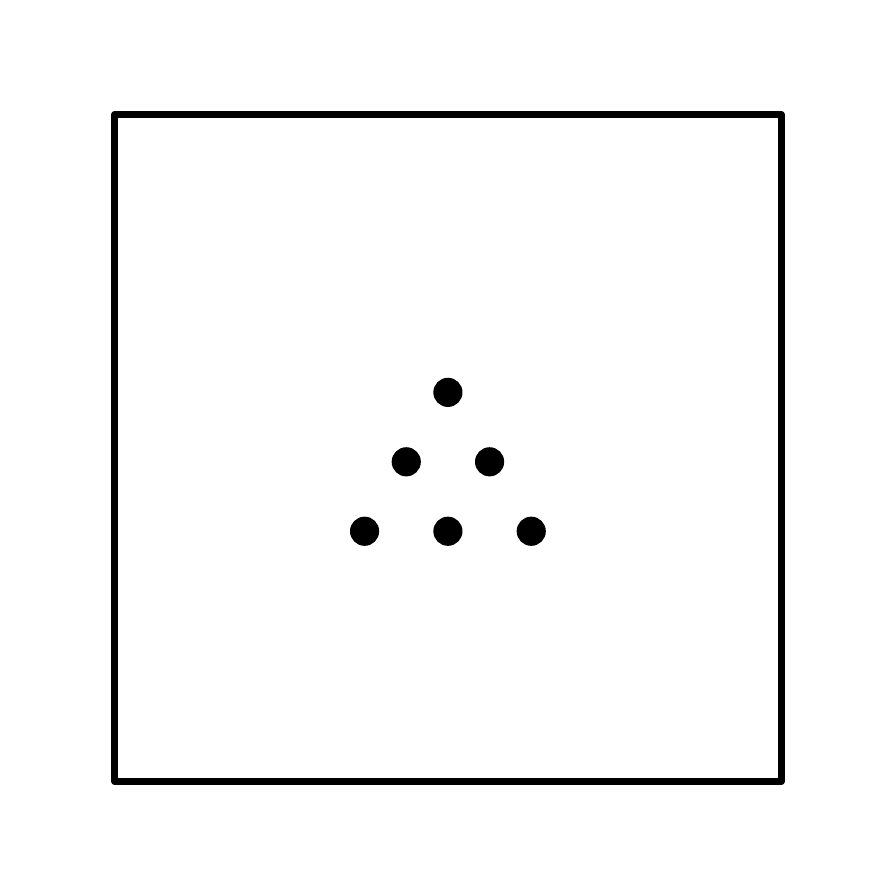}} \\
\dim\calS_3^-\Lambda^0(\square_2)=12 & \dim\calS_3^-\Lambda^1(\square_2)=17 & \dim\calS_3^-\Lambda^2(\square_2)=6 \\
}
\end{equation*}
\caption{Element diagrams for $\calS_2^-\Lambda^k(\square_2)$ and $\calS_3^-\Lambda^k(\square_2)$, shown as subcomplexes of the deRham complex for $\square_2$.
Each dot or arrow in the diagrams indicates a degree of freedom associated to that piece of the geometry (vertex, edge, or interior).
}
\label{fig:Sm2and3k2}
\end{figure}

The $\calS_r^-\Lambda^k(\square_n)$ elements of most immediate relevance to applications are those for small values of $n$ and $r$.
We now examine some of these cases in greater detail, using a mix of exterior calculus and vector calculus notation.
Formal definitions of the notation and generalized formulae using exclusively exterior calculus notation are given in Sections~\ref{sec:bkgd}-\ref{sec:unisolv-min} and a description of how to convert between vector and exterior calculus notation is given in Appendix~\ref{appdx:vec-calc}.

\textbf{The spaces $\calS_r^-\Lambda^k(\square_2)$}.
The element diagrams in Figure~\ref{fig:Sm2and3k2} indicate the association of degrees of freedom to portions of the square geometry for $\calS_2^-\Lambda^k(\square_2)$ and $\calS_3^-\Lambda^k(\square_2)$.
The degrees of freedom for $\calS_r^-\Lambda^1(\square_2)$ are 
\begin{align*}
u\longmapsto \int_e u\cdot \vec t~p, &\qquad p\in\calP_{r-1}(e),~~\text{$e$ an edge of $\square_2$ with unit tangent $\vec t$},\\[2mm]
u\longmapsto \int_{\square_2} u\cdot \vec p, &\qquad \vec p\in\left[\calP_{r-3}(\square_2)\right]^2\oplus\grad\calH_{r-1}\Lambda^0(\square_2).
\end{align*}
The notation $\grad\calH_{r-1}\Lambda^0(\square_2)$ above should be interpreted as the vector proxies for the exterior derivative applied to homogenous polynomials of degree $r-1$ in two variables.
Observe that if we exclude only the degrees of freedom associated to $\grad\calH_{r-1}\Lambda^0(\square_2)$, we are left with the degrees of freedom for the regular serendipity space $\calS_{r-1}\Lambda^1(\square_2)$.\\

\begin{figure}[t]
\centering
\begin{equation*}
\xymatrix @R=.02in{
\parbox{.18\textwidth}
{\includegraphics[width=.18\textwidth]{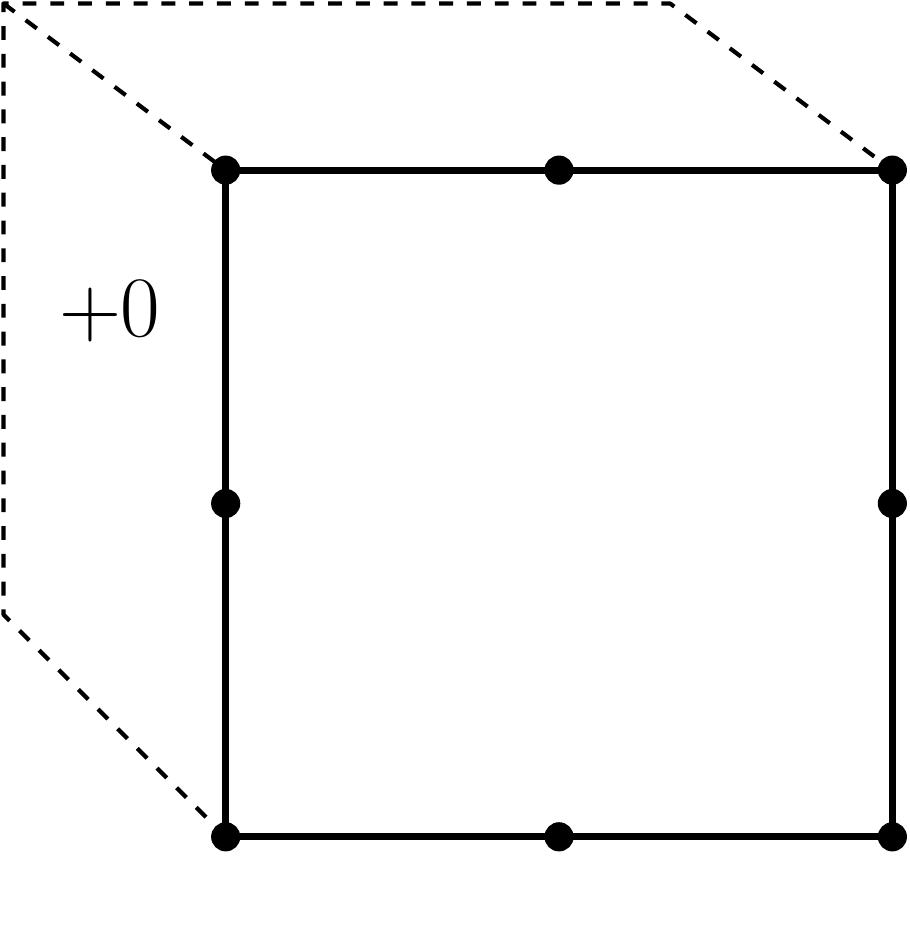}}  \ar[r]^-{d_0}_{\grad}
&
\parbox{.18\textwidth}
{\includegraphics[width=.18\textwidth]{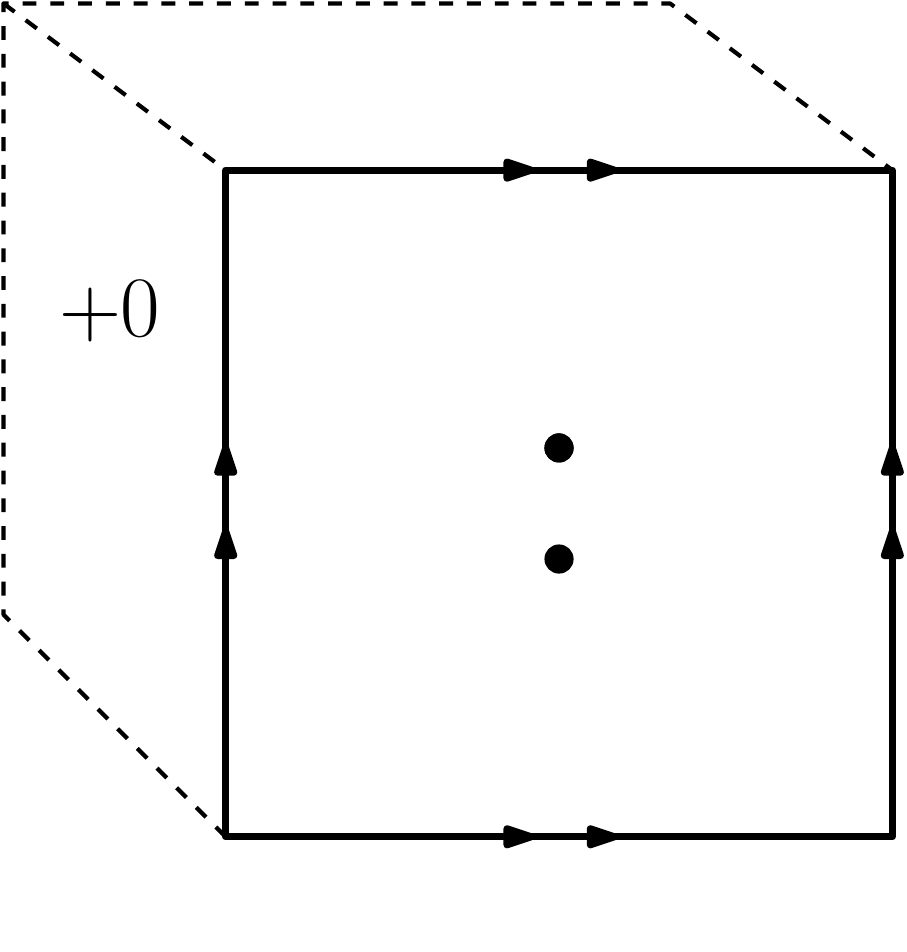}} \ar[r]^-{d_1}_{\curl}
&
\parbox{.18\textwidth}
{\includegraphics[width=.18\textwidth]{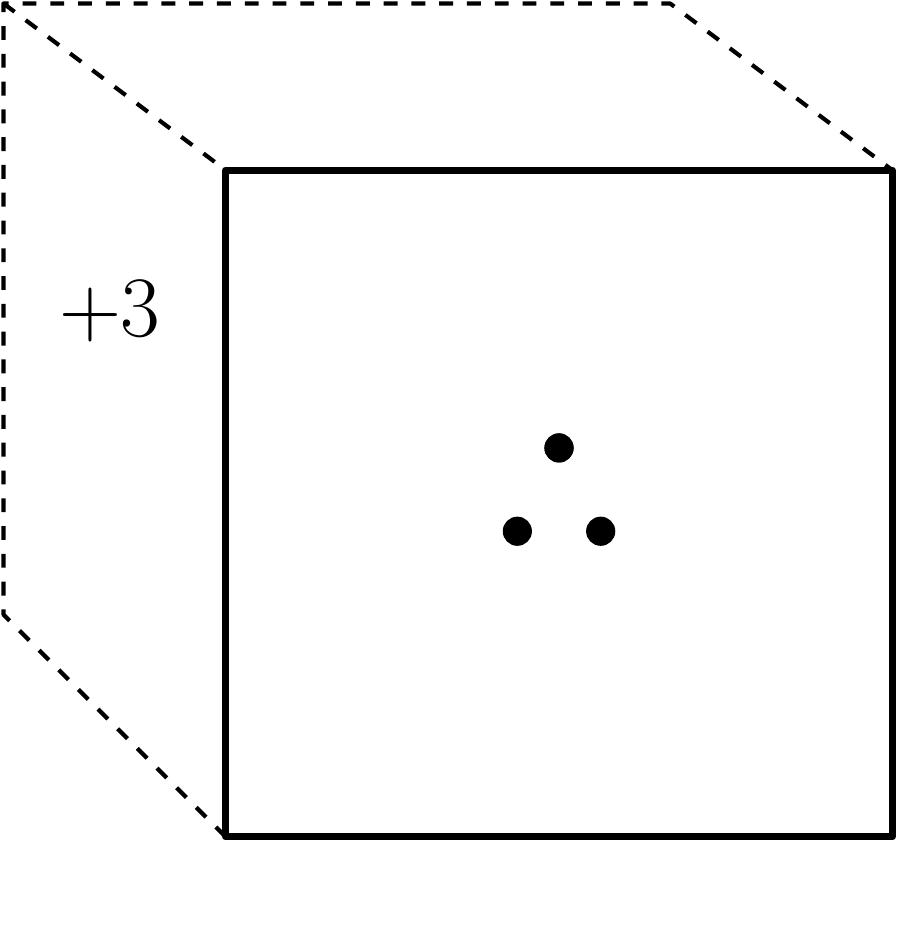}} \ar[r]^-{d_2}_{\div}
&
\parbox{.18\textwidth}
{\includegraphics[width=.18\textwidth]{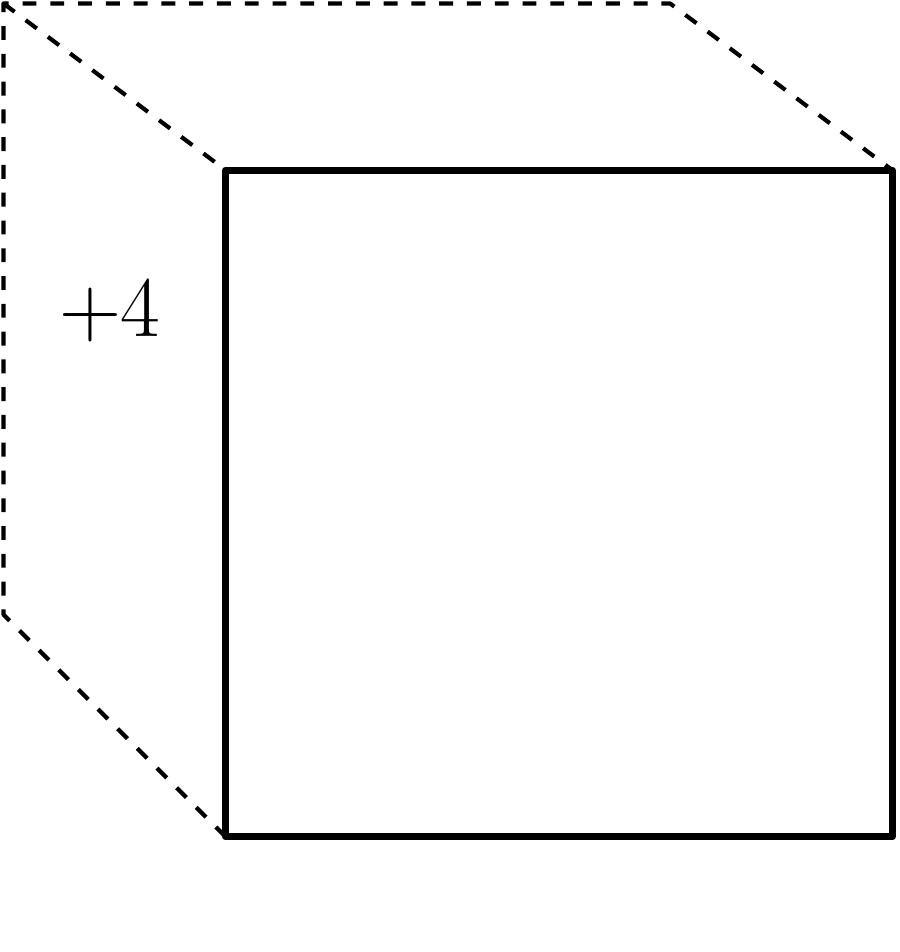}}\\
\calS_2^-\Lambda^0(\square_3) & \calS_2^-\Lambda^1(\square_3) & \calS_2^-\Lambda^2(\square_3) & \calS_2^-\Lambda^3(\square_3)\\
\dim = 20 & \dim =36 & \dim =21 & \dim=4
}
\end{equation*}
\caption{$\calS_2^-\Lambda^k(\square_2)$ for $k=0,1,2,3$.
Degrees of freedom on a representative face are shown, along with a count of $+X$ for the number of degrees of freedom associated to the interior of the cube.
}
\label{fig:Sm2k3}
\end{figure}

\textbf{The spaces $\calS_r^-\Lambda^k(\square_3)$}.
Moving to cubes, element diagrams for the $\calS_2^-\Lambda^k(\square_3)$ spaces are shown in Figure~\ref{fig:Sm2k3}.
In these figures, degrees of freedom associated to vertices, edges, or faces of the cube are shown on the front face only while the number of degrees of freedom associated to the interior of the cube are indicated by $+X$.
Looking only at the front face degrees of freedom in Figure~\ref{fig:Sm2k3} for $k=0,1,2$, we see exactly the same sequence as shown in the top row of Figure~\ref{fig:Sm2and3k2}, reflecting the fact that the $\calS_r^-\Lambda^k(\square_n)$ spaces have the trace property.
We also observe that $\calS_2^-\Lambda^0(\square_3)=\calS_2\Lambda^0(\square_3)$ and $\calS_2^-\Lambda^3(\square_3)=\calS_1\Lambda^3(\square_3)$.
Further, the lowest order spaces also coincide with the tensor product differential form spaces, i.e. $\calS_1^-\Lambda^k(\square_3)=\calQ_1^-\Lambda^k(\square_3)$ for $k=0,1,2,3$.

The degrees of freedom for $\calS_r^-\Lambda^1(\square_3)$ are 
\begin{align*}
u\longmapsto \int_e u\cdot \vec t~p, &\qquad p\in\calP_{r-1}(e),~~\text{$e$ an edge of $\square_3$ with unit tangent $\vec t$},\\[2mm]
u\longmapsto \int_{f} (u\times \hat n) \cdot \vec p, &\qquad \vec p\in\left[\calP_{r-3}(f)\right]^2\oplus\grad\calH_{r-1}\Lambda^0(f),\\
 &\qquad\qquad\qquad \text{$f$ a face of $\square_3$ with unit normal $\hat n$},\\[2mm]
u\longmapsto \int_{\square_3} u\cdot \vec p, &\qquad \vec p\in\left[\calP_{r-5}(\square_3)\right]^3\oplus\curl\calH_{r-3}\Lambda^1(\square_3).
\end{align*}
As in the $n=2$ case, we observe that removing the degrees of freedom associated to $\grad\calH_{r-1}\Lambda^0(f)$ and $\curl\calH_{r-3}\Lambda^1(\square_3)$ leaves only the degrees of freedom for $\calS_{r-1}\Lambda^1(\square_3)$.
 
The degrees of freedom for $\calS_r^-\Lambda^2(\square_3)$ are 
\begin{align*}
u\longmapsto \int_{f} u\cdot \hat n\ p, &\qquad p\in \calP_{r-1}(f),~\text{$f$ a face of $\square_3$ with unit normal $\hat n$},\\[2mm]
u\longmapsto \int_{\square_3} u\cdot p, &\qquad \vec p\in\left[\calP_{r-3}(\square_3)\right]^3\oplus\grad\calH_{r-1}\Lambda^0(\square_3).
\end{align*}
Again, excluding the degrees of freedom associated to $\grad\calH_{r-1}\Lambda^0(\square_3)$, we are left with the degrees of freedom for $\calS_{r-1}\Lambda^2(\square_3)$.\\

The remainder of the paper is organized as follows.
In Section~\ref{sec:bkgd}, we review relevant background and notation from finite element exterior calculus and compare the trimmed serendipity elements to other elements in the literature.
In Section~\ref{sec:srm-def}, we prove various properties of the $\calS_r^-\Lambda^k(\bbR^n)$ spaces, including a formula to compute their dimension.
In Section~\ref{sec:unisolv-min}, we state a set of degrees of freedom and prove they are unisolvent for $\calS_r^-\Lambda^k(\square_n)$.
We also explain and establish minimality in the context of compatible finite element systems.
Finally, we summarize the key results of our work and give an outlook on the future directions they suggest in Section~\ref{sec:conc}.
Appendix~\ref{appdx:vec-calc} provides a detailed description of the relation between exterior calculus and vector calculus notation in the context studied here.

\section{Notation and Relation to Prior Work}
\label{sec:bkgd}

We use the same notation as Arnold and Awanou~\cite{AA2014} and will now review the relevant definitions to aid in comparison to prior work.
Let $\alpha\in\bbN^n$ be a multi-index and let $\sigma$ be a subset of $\{1,\ldots,n\}$ consisting of $k$ distinct elements $\sigma(1),\ldots,\sigma(k)$ with $0\leq k\leq n$.
The \textit{form monomial} $x^\alpha dx_\sigma$ is the differential $k$-form on $\bbR^n$ given by
\begin{equation}
\label{eq:form-mon}
x^{\alpha}dx_{\sigma}:=\left(x_1^{\alpha_1}x_2^{\alpha_2}\dots x_n^{\alpha_n}\right)dx_{\sigma(1)}\wedge\dots\wedge dx_{\sigma(k)}.
\end{equation}
The \textit{degree} of $x^{\alpha}dx_{\sigma}$ is $|\alpha|:=\sum_{i=1}^n\alpha_i$.
The space of differential $k$-forms with polynomial coefficients of homogeneous degree $r$ is denoted $\calH_r\Lambda^k(\bbR^n)$.
A basis for this space is the set of form monomials such that $|\alpha|=r$ and $|\sigma|=k$.
The exterior derivative $d$ and Koszul operator $\kappa$ are maps
\[d:\calH_r\Lambda^k(\bbR^n)\raw\calH_{r-1}\Lambda^{k+1}(\bbR^n)\quad \kappa:\calH_r\Lambda^k(\bbR^n)\raw\calH_{r+1}\Lambda^{k-1}(\bbR^n).\]
In coordinates, they are defined on form monomials by
\begin{align}
d(x^\alpha dx_\sigma) & := \sum_{i=1}^n \left(\frac{\partial x^\alpha}{\partial x_i}dx_i\right)\wedge dx_{\sigma(1)}\wedge\dots\wedge dx_{\sigma(k)}, \label{eq:def-d}\\
\kappa(x^\alpha dx_\sigma) & := \sum_{i=1}^k \left((-1)^{i+1}x^\alpha x_{\sigma(i)}\right)dx_{\sigma(1)}\wedge\cdots\wedge\widehat{dx_{\sigma(i)}}\wedge\cdots\wedge dx_{\sigma(k)}.\label{eq:def-kappa}
\end{align}
The notation $\widehat{dx_{\sigma(i)}}$ indicates that the term is omitted from the wedge product.
We will make frequent use of the \textit{homotopy formula} in this context~\cite[Theorem 3.1]{AFW2006}, which is also called \textit{Cartan's magic formula}:
\begin{equation}
\label{eq:magic}
(d\kappa+\kappa d)\omega = (r+k)\omega,\quad\omega\in\calH_r\Lambda^k(\bbR^n).
\end{equation}
As shown in~\cite[equation (3.10)]{AFW2006}, it follows that 
\begin{equation}
\label{eq:hr-decomp}
\calH_r\Lambda^k(\bbR^n) = \kappa\calH_{r-1}\Lambda^{k+1}(\bbR^n)\oplus d\calH_{r+1}\Lambda^{k-1}(\bbR^n).
\end{equation}
The space of polynomial differential $k$-forms of degree at most $r$ is
\begin{equation}
\label{eq:pr-sum-of-hr}
\calP_r\Lambda^k(\bbR^n):=\bigoplus_{j=0}^r\calH_j\Lambda^k(\bbR^n).
\end{equation}
The definitions of $d$ and $\kappa$ extend linearly over $\calP_r\Lambda^k$.
As a consequence, 
\begin{equation}
\label{eq:pr-decomp}
\calP_r\Lambda^k(\bbR^n) = \kappa\calP_{r-1}\Lambda^{k+1}(\bbR^n)\oplus d\calP_{r+1}\Lambda^{k-1}(\bbR^n).
\end{equation}
Observe that if $\omega\in\calP_r\Lambda^k(\bbR^n)$ can be written as both an image of $\kappa$ and as an image of $d$, then $\omega=0$.

The ``trimmed'' space of polynomial differential $k$-forms of degree at most $r$ is
\begin{equation}
\label{eq:prm-def}
\calP_r^-\Lambda^k(\bbR^n) := \calP_{r-1}\Lambda^k(\bbR^n) \oplus \kappa\calH_{r-1}\Lambda^{k+1}(\bbR^n).
\end{equation}
The relation of the $\calP_r\Lambda^k(\bbR^n)$ and $\calP_r^-\Lambda^k(\bbR^n)$ spaces to the well-known \nedelec~\cite{N1980,N1986}, Raviart-Thomas~\cite{RT1977} and Brezzi-Douglas-Marini~\cite{BDM85} elements on simplices is described in the work of Arnold, Falk and Winther~\cite{AFW2006,AFW2010} and summarized in the Periodic Table of the Finite Elements~\cite{AL2014}.

An essential precursor to the finite element exterior calculus framework just described is the work of Hiptmair~\cite{Hip99}, in which the spaces $\calP_r^-\Lambda^k$ were introduced under different notation.
In place of the Koszul operator, Hiptmair uses a potential mapping, $k_{\bfa}$, which satisfies the formula $d(k_\bfa(\omega))+k_\bfa(d\omega)=\omega$, similar to (\ref{eq:magic}) but without the factor of $(r+k)$.
The mapping $k_\bfa$ is defined in terms of a point $\bfa$ inside a star-shaped domain whereas the Koszul operator implicitly chooses the origin as a reference point.
Since there is no true inverse for the $d$ operator, using $k_\bfa$ in place of $\kappa$ could provide some additional insight, however, we have found the $\kappa$ operator to be quite natural for characterizing the structure of polynomial finite element differential form spaces.

To build the serendipity spaces on $n$-dimensional cubes, we need some additional definitions.
Let $\sigma^\ast$ denote the complement of $\sigma$, i.e. $\sigma^\ast:=\{1,\ldots,n\}-\sigma$.
The \textit{linear degree} of $x^{\alpha}dx_{\sigma}$ is defined to be
\begin{equation}
\label{eq:lindeg-def}
\ldeg(x^\alpha dx_\sigma):=\#\{i\in\sigma^\ast~:~\alpha_i=1\}.
\end{equation}
Put differently, the linear degree of $x^\alpha dx_\sigma$ counts the number of entries in $\alpha$ equal to 1, excluding entries whose indices appear in $\sigma$.
Note that if $k=0$ then $\sigma=\varnothing$ and there is no `exclusion' in the counting of linear degree.
Likewise, if $k=n$ then $\sigma^\ast=\varnothing$ and $\ldeg(x^\alpha dx_\sigma)=0$ for any $\alpha$.
The linear degree of the sum of two or more form monomials is defined as the minimum of the linear degrees of the summands.

The subset of $\calH_r\Lambda^k(\bbR^n)$ that has linear degree at least $\ell$ is denoted
\begin{equation}
\label{eq:hrlk-def}
		\calH_{r,l}\Lambda^k(\bbR^n):=\left\{\omega\in\calH_r\Lambda^k(\bbR^n)~|~\mbox{ldeg }\omega\geq l\right\}.
\end{equation}
A key building block for both the serendipity and trimmed serendipity spaces is
\begin{equation}
\label{eq:jrk-def}
\calJ_r \Lambda^k(\bbR^n) := \sum_{l \geq 1} \kappa \calH_{r+l-1, l} \Lambda^{k+1}(\bbR^n).
\end{equation}
The following proposition gives a simple and useful characterization of $\calJ_r \Lambda^k(\bbR^n)$.

\begin{prop}[{\cite[Proposition 3.1]{AA2014}}]
\label{prop:jrk-char}
The space $\calJ_r^k(\bbR^n)$ is the span of $\kappa m$ for all $(k+1)$--form monomials with $\deg m\geq r$ and $\deg m-\ldeg~m\leq r-1$.
\end{prop}

Note that every element in $\calJ_r\Lambda^k(\bbR^n)$ lies in the range of $\kappa$.
Using this fact, we develop some basic results about $\calJ_r\Lambda^k(\bbR^n)$ that will be useful in our development of the $\calS_r^-\Lambda^k$ spaces.
In the proof of~\cite[Theorem 3.4]{AA2014}, it is shown that 
\begin{equation}
\label{eq:jr-in-pandj}
\calJ_r\Lambda^k(\bbR^n)\subset\calP_{r+1}\Lambda^k(\bbR^n) + \calJ_{r+1}\Lambda^k(\bbR^n).
\end{equation}
We can exclude images by $d$ from the right side, yielding
\begin{equation}
\label{eq:jr-in-pandj-sharp}
\calJ_r\Lambda^k(\bbR^n)\subset\kappa\calP_{r}\Lambda^{k+1}(\bbR^n) + \calJ_{r+1}\Lambda^k(\bbR^n).
\end{equation}
Further, by (\ref{eq:magic}), we have that $(d\kappa+\kappa d)\calJ_r\Lambda^k(\bbR^n)=\calJ_r\Lambda^k(\bbR^n)$.
Since $\kappa\kappa=0$, we have $d\kappa\calJ_r\Lambda^k(\bbR^n)=0$, and thus
\begin{equation}
\label{eq:kdJ-J}
\kappa d\calJ_r\Lambda^k(\bbR^n) = \calJ_r\Lambda^k(\bbR^n).
\end{equation}
\indent
The space of serendipity differential $k$-forms of order $r$ is given by
\begin{equation}
\label{eq:srk-decomp}
\calS_r\Lambda^k(\bbR^n)= \calP_r\Lambda^k(\bbR^n) \oplus \calJ_r \Lambda^k(\bbR^n)\oplus  d \calJ_{r+1}\Lambda^{k-1}(\bbR^n).
\end{equation}
The fact that this sum is direct is proven in~\cite{AA2014}.
Note that the second summand vanishes when $k=n$, since $\Lambda^{n+1}(\bbR^n)=0$ while the third summand vanishes when $k=0$, since $\Lambda^{-1}=0$ by definition.
Given $x^\alpha dx_\sigma\in\calS_r\Lambda^k(\bbR^n)$, the degree property from~\cite[Theorem 3.2]{AA2014} ensures that
\begin{equation}
\label{eq:sr-deg-prop}
\deg(x^\alpha dx_\sigma)\leq r+n-k-\delta_{k0}\quad\text{and}\quad\deg(x^\alpha dx_\sigma)-\ldeg(x^\alpha dx_\sigma)\leq r+1-\delta_{k0},
\end{equation}
where $\delta_{ij}$ denotes the Kronecker delta function.

The serendipity spaces satisfy an inclusion property~\cite[Theorem 3.4]{AA2014}:
\begin{equation}
\label{eq:sr-inc}
\calS_{r}\Lambda^k(\bbR^n) \subset  \calS_{r+1}\Lambda^k(\bbR^n),
\end{equation}
and a subcomplex property~\cite[Theorem 3.3]{AA2014}:
\begin{equation}
\label{eq:sr-subcpx}
d\calS_{r+1}\Lambda^{k-1}(\bbR^n)\subset\calS_r\Lambda^k(\bbR^n).
\end{equation}
They also satisfy a containment property with respect to $\kappa$, namely,
\begin{equation}
\label{in:srk}
\kappa \calS_{r-1} \Lambda^{k}(\bbR^n) \subset \calS_r \Lambda^{k-1}(\bbR^n). 
\end{equation}
The proof is a direct consequence of (\ref{eq:pr-decomp}), (\ref{eq:kdJ-J}), and (\ref{eq:srk-decomp}).

The $\calS_r\Lambda^k(\bbR^n)$ spaces can be collected into a cochain complex with decreasing $r$, which we denote by $\calS_{r-\bs}\Lambda^\bs$.
The resulting sequence, as well as those for $\calP_{r-\bs}\Lambda^\bs$ and $\calP_{r}^-\Lambda^\bs$, augmented by $\bbR$ in front of the first term, are all exact.
Written out, these are
\begin{alignat}{6}
&0 \to \bbR\to\calS_{r} \Lambda^0 &\to~ & \calS_{r-1} \Lambda^{1} &\to & ~\cdots~ \to \calS_{r-n+1} \Lambda^{n-1} &&\to \calS_{r-n} \Lambda^{n} &&\to 0, \label{eq:srk-seq}\\[1mm]
& 0 \to \bbR\to\calP_{r} \Lambda^0 &\to~ & \calP_{r-1} \Lambda^{1} &\to & ~\cdots~ \to \calP_{r-n+1} \Lambda^{n-1} &&\to  \calP_{r-n} \Lambda^{n} &&\to 0, \\[1mm]
& 0 \to \bbR\to\calP_{r}^- \Lambda^0 &\to ~& \calP_{r}^- \Lambda^{1} &\to & ~\cdots~ \to \calP_{r}^- \Lambda^{n-1} &&\to  \calP_{r}^- \Lambda^{n} &&\to 0. \label{eq:prm-seq}
\end{alignat}
All the above sequences can serve as finite element subcomplexes of the deRham complex for a domain.
Such subcomplexes help guide the selection of pairs of spaces for mixed finite element methods that have guaranteed stability and convergence properties.

\subsection*{Comparison to prior and contemporary work.}
As mentioned in the introduction, there has been a recent spate of research into conforming finite elements on meshes of $n$-dimensional cubes.
The $\calS_r\Lambda^k(\bbR^n)$ family was defined first in the $H^1$-conforming ($k=0$) case in~\cite{AA2011} and then for any $0\leq k\leq n$ in~\cite{AA2014}.
The relation of the $\calS_r\Lambda^1(\bbR^2)$ elements to the Brezzi-Douglas-Marini~\cite{BDM85} elements on rectangles is described in~\cite{AA2014} and by the Periodic Table of the Finite Elements~\cite{AL2014}.
The relation between the trimmed and non-trimmed serendipity families is described by Lemma~\ref{lem:sr-alt-decomp} below.
For $0<k\leq n$, we will see that $\dim\calS_r^-\Lambda^k(\bbR_n)<\dim\calS_r\Lambda^k(\bbR_n)$, indicating that the trimmed and non-trimmed serendipity families are truly distinct.

By converting from exterior caclulus to vector calculus notation, we can  identify the relation of the $\calS_r^-\Lambda^k(\bbR^n)$ spaces to finite element spaces defined in recent work by Arbogast and Correa~\cite{AC2015} and by Cockburn and Fu~\cite{CF2016}.
Both works examine various families of elements, and each work presents one family that is essentially the same as the trimmed serendipity elements, as explained in the following Propositions.
Note that both sets of authors use $k$ to indicate polynomial degree, but we have changed the notation to $r$ to match the conventions of finite element exterior calculus.
We use the notation $\square_n$ to denote the cube $[-1,1]^n\subset\bbR^n$.

\begin{prop}
\label{prop:ac-comp}
Define the pair of spaces $\left(\mathbf{V}_{AC}^r,W_{AC}^r\right)\subset H(\div,\square_2)\times L^2(\square_2)$ as in~\cite{AC2015}.
Let $\rot\mathbf{V}_{AC}^r$ denote the application of the $\rot$ operator to all the vectors in $\mathbf{V}_{AC}^r$, which has the effect of rotating each vector in the field by $\pi/2$.
Then, interpreted as differential forms via the flat operator, $\left(\rot\mathbf{V}_{AC}^r,W_{AC}^r\right)$ is identical to $\left(\calS_{r+1}^-\Lambda^1(\square_2),\calS_{r+1}^-\Lambda^2(\square_2)\right)$.
\end{prop}

\begin{prop}
\label{prop:cf-comp}
The sequence of spaces denoted $S_{2,r}^{\square_2}(K)$ in~\cite[Theorem 3.3]{CF2016}, interpreted as differential forms via the flat operator, is identical to the sequence
\[\calS_{r+1}^-\Lambda^0(\square_2) \to \calS_{r+1}^-\Lambda^1(\square_2) \to \calS_{r+1}^-\Lambda^2(\square_2).\]
Further, the sequence denoted $S^{\square_3}_{2,r}$ in~\cite[Theorem 3.6]{CF2016}, is identical to the sequence
\[\calS_{r+1}^-\Lambda^0(\square_3) \to \calS_{r+1}^-\Lambda^1(\square_3) \to \calS_{r+1}^-\Lambda^2(\square_3) \to \calS_{r+1}^-\Lambda^3(\square_3).\]
\end{prop}

Detailed proofs of both propositions are given in Appendix~\ref{appdx:vec-calc}.
We can now make a precise statement about the novelty of the trimmed serendipity spaces.
The spaces $\calS_r^-\Lambda^k(\square_n)$ can be recognized as differential form analogues of ($i$) the mixed finite element method presented in~\cite{AC2015} when applied to affinely-mapped square meshes (as opposed to general quadrilateral meshes), and ($ii$) the second of the four families of elements on squares and cubes presented in~\cite{CF2016}.
For $n\geq 4$, the trimmed serendipity spaces are entirely new to the literature, modulo the fact that the $k=0$ and $k=n$ cases reduce to the non-trimmed serendipity spaces as described in Lemma~\ref{lem:sr-alt-decomp}.

Further comparison can be made in regards to degrees of freedom.
We define degrees of freedom for $\calS_r^-\Lambda^k(\square_n)$ in equation (\ref{eq:Srm-dofs}) and prove they are unisolvent for $\calS_r^-\Lambda^k(\square_n)$ in Theorem~\ref{thm:unisolv}.
The degrees of freedom given by Arbogast and Correa~\cite{AC2015} are slightly different, in that they are indexed in part by vectors of polynomials that vanish on certain edges of $\square_2$, whereas our degrees of freedom are indexed by spaces of polynomial differential forms without regard to the basis used to define them.
The spaces of Cockburn and Fu~\cite{CF2016} are not equipped with degrees of freedom and so no comparison is possible in this case.
Finally, we mention the virtual element space $VEMS^f_{r,r,r-1}$, recently defined in work by Beir\~{a}o da Veiga et al~\cite{BeiEtAl16}.
The number of degrees of freedom for this space appears to agree with the number of degrees of freedom for $\calS_{r+1}^-\Lambda^1(\square_2)$ in the case of a square, the main difference being a vector calculus treatment of indexing spaces in place of the differential form terminology used here.
Since the virtual element method does not employ spaces of local basis functions, further comparison between the methods is a larger question for future work.

\section{The $\calS_r^-\Lambda^k$ spaces}
\label{sec:srm-def}

We define the trimmed serendipity spaces for $r \geq 1$, $k \geq 0$ by
\begin{equation}
\label{eq:srmk-def}
\calS_r^-\Lambda^k(\bbR^n) :=\calS_{r-1}\Lambda^k(\bbR^n) +\kappa\calS_{r-1}\Lambda^{k+1}(\bbR^n).
\end{equation}
The trimmed serendipity spaces share many analogues with the trimmed polynomial spaces, as we now establish.
Throughout, we fix the top dimension to be $n\geq 1$ and omit the notation $(\bbR^n)$, except when it is needed for clarity.

\begin{theorem}[Inclusion property]
Let $n,r\geq 1$, and $0 \leq k \leq n$.
Then
\begin{equation}
\label{eq:inclusion}
\calS_{r}\Lambda^k \subset \calS_{r+1}^-\Lambda^k \subset \calS_{r+1}\Lambda^k.
\end{equation}
\end{theorem}
\begin{proof}
The first inclusion is immediate from (\ref{eq:srmk-def}).
For the second inclusion, the inclusion property (\ref{eq:sr-inc}) implies that $\calS_{r}\Lambda^k \subset  \calS_{r+1}\Lambda^k$.
Hence we only need to show that $\kappa\calS_{r}\Lambda^{k+1}\subset\calS_{r+1}\Lambda^k$.
Decomposing $\kappa\calS_{r}\Lambda^{k+1}$ by (\ref{eq:srk-decomp}), we have $\kappa\calP_{r}\Lambda^{k+1}\subset\calP_{r+1}\Lambda^k\subset\calS_{r+1}\Lambda^k$, $\kappa \calJ_{r} \Lambda^{k+1}=0$, and, by (\ref{eq:kdJ-J}), $\kappa d\calJ_{r+1}\Lambda^k= \calJ_{r+1}\Lambda^k\subset\calS_{r+1}\Lambda^k$, thus completing the proof.
\end{proof}

\begin{theorem}[Subcomplex property]
\label{thm:subcpx}
Let $n,r \geq 1$, and $0 < k \leq n$.
Then
\begin{equation}
\label{eq:subcpx}
d\calS_r^-\Lambda^k \subset \calS_r^-\Lambda^{k+1}.
\end{equation}
\end{theorem}
\begin{proof}
Using (\ref{eq:inclusion}) and (\ref{eq:sr-subcpx}), we have $d\calS_r^-\Lambda^k \subset d\calS_r\Lambda^k \subset \calS_{r-1}\Lambda^{k+1}\subset \calS_r^-\Lambda^{k+1}$.
\end{proof}

We now develop a direct sum decomposition of $\calS_r^-\Lambda^k$ whose proof is straightforward by virtue of an analogous decomposition of $\calS_r\Lambda^k$.

\begin{theorem}[Direct sum decomposition]
\label{thm:srm-directsumdecomp}
Let $n,r \geq 1$, and $0 \leq k \leq n$. 
Then $\calS_r^-\Lambda^k$, as defined by (\ref{eq:srmk-def}), can also by written as the direct sum
\begin{equation}
\label{eq:srm-directsumdecomp}
\calS_r^-\Lambda^k = \calP_r^-\Lambda^k \oplus \calJ_r\Lambda^k \oplus d \calJ_{r} \Lambda^{k-1}.
\end{equation}
Further, any element $\omega\in\calS_r^-\Lambda^k$ can be written as $\omega=d\alpha+\kappa \beta$ where $d\alpha\in\calS_{r-1}\Lambda^k$ and $\kappa\beta\in\calS_r\Lambda^k$.
In particular, $\alpha\in \calP_r\Lambda^{k-1}\oplus \calJ_{r} \Lambda^{k-1}$ and $\beta\in \calP_{r-1}\Lambda^{k+1} \oplus \sum_{l \geq 1}  \calH_{r+l-1, l} \Lambda^{k+1}$.
\end{theorem}

\begin{proof}
First expand (\ref{eq:srmk-def}) via (\ref{eq:srk-decomp}).
Since $\kappa^2 = 0$, we have
\[
\calS_r^-\Lambda^k = (\calP_{r-1}\Lambda^k  \oplus \calJ_{r-1} \Lambda^k  \oplus  d \calJ_{r} \Lambda^{k-1}) + (\kappa\calP_{r-1}\Lambda^{k+1}  \oplus\kappa d\calJ_r\Lambda^k). 
\]
By (\ref{eq:kdJ-J}), we can replace $\kappa d\calJ_r\Lambda^k$ by $\calJ_r\Lambda^k$.
Further, by (\ref{eq:jr-in-pandj-sharp}), applied to the $\calJ_{r-1}\Lambda^k$ term and re-ordering, we now have
\[
\calS_r^-\Lambda^k = \calP_{r-1}\Lambda^k + \kappa\calP_{r-1}\Lambda^{k+1} + \calJ_r\Lambda^k + d \calJ_{r} \Lambda^{k-1} .
\]
The first two terms summands give $\calP_{r-1}\Lambda^k + \kappa\calP_{r-1}\Lambda^{k+1}=\calP_r^-\Lambda^k$, which establishes (\ref{eq:srm-directsumdecomp}) as a summation formula.

We now show that (\ref{eq:srm-directsumdecomp}) is direct.
Observe that $\calP_r^-\Lambda^k=\kappa\calP_{r-1}\Lambda^{k+1} \oplus d\calP_r\Lambda^{k-1}$, with directness of the sum following from the 
observation after (\ref{eq:pr-decomp}) that 0 is the only polynomial differential form in the image of both $\kappa$ and $d$.
Therefore, the intersection of $\kappa\calP_{r-1}\Lambda^{k+1} + \calJ_r\Lambda^k$ with  $d\calP_{r}\Lambda^k +  d \calJ_{r} \Lambda^{k-1}$ is $\{0\}$.
Now, elements of $\kappa\calP_{r-1}\Lambda^{k+1}$ are of degree at most $r$ while elements of $\calJ_r\Lambda^k$ are of degree at least $r+1$.
Similarly, elements of $d\calP_{r}\Lambda^k$ are of degree at most $r-1$ while elements of $d \calJ_{r} \Lambda^{k-1}$ are of degree at least $r$.
Hence both pairs are direct sums and (\ref{eq:srm-directsumdecomp}) is established.

For the last statement, again consider the direct sum $\calP_r^-\Lambda^k= \kappa\calP_{r-1}\Lambda^{k+1}\oplus d\calP_r\Lambda^{k-1}$.
Given $\omega\in\calS_r^-\Lambda^k$, we can thus write $\omega=d\alpha+\kappa\beta$ such that $d\alpha\in d\calP_r\Lambda^{k-1}\oplus d\calJ_{r} \Lambda^{k-1}$ and $\kappa\beta\in\kappa\calP_{r-1}\Lambda^{k+1}\oplus\calJ_r\Lambda^k$.
We have $d\calP_r\Lambda^{k-1}\oplus d\calJ_{r} \Lambda^{k-1}\subset\calP_{r-1}\Lambda^k\oplus d\calJ_{r} \Lambda^{k-1}\subset\calS_{r-1}\Lambda^k$ and $\kappa\calP_{r-1}\Lambda^{k+1} \oplus\calJ_r\Lambda^k\subset \calP_r\Lambda^k\oplus\calJ_r\Lambda^k\subset\calS_r\Lambda^k$, as seen from (\ref{eq:srk-decomp}).
\end{proof}

\begin{lemma}
\label{lem:sr-alt-decomp}
Let $n,r\geq 1$.
\begin{enumerate}
\renewcommand{\theenumi}{\roman{enumi}}
\item $\calS_r^-\Lambda^0=\calS_r\Lambda^0$, \\[-3mm]
\item $\calS_r^-\Lambda^n=\calS_{r-1}\Lambda^n$, \\[-3mm]
\item $\calS_r^-\Lambda^k + d\calS_{r+1}\Lambda^{k-1}=\calS_r\Lambda^k$.
\end{enumerate}
\end{lemma}

\begin{proof}
For (i), note that $\calP_r\Lambda^0  = \kappa\calP_{r-1}\Lambda^1$ by (\ref{eq:pr-decomp}) and $\kappa d\calJ_r\Lambda^0 = \calJ_r\Lambda^0$ by (\ref{eq:kdJ-J}).
We decompose $\calS_r^-\Lambda^k$ according to (\ref{eq:srmk-def}) and then decompose the summands according to (\ref{eq:srk-decomp}), yielding
\begin{align*}
\calS_r^-\Lambda^0 
 & = \calS_{r-1}\Lambda^0 + \kappa\calS_{r-1}\Lambda^1 \\
 & = \left(\calP_{r-1}\Lambda^0 + \calJ_{r-1}\Lambda^0\right) + \left(\kappa\calP_{r-1}\Lambda^1 + \kappa d\calJ_r\Lambda^0\right) \\
 & = \calP_r\Lambda^0 + \calJ_{r-1}\Lambda^0 +\calJ_r\Lambda^0.
\end{align*}
By (\ref{eq:jr-in-pandj}), $\calJ_{r-1}\Lambda^0\subset\calP_r\Lambda^0+\calJ_r\Lambda^0$ and so $\calS_r^-\Lambda^0(\bbR^n) = \calP_r\Lambda^0 +  \calJ_r\Lambda^0 = \calS_r\Lambda^0$.
Part (ii) is an immediate consequence of (\ref{eq:srmk-def}), since there are no $(n+1)$-forms on $\bbR^n$.

For (iii), we have $\calS_r^-\Lambda^k \subset\calS_r\Lambda^k$ by (\ref{eq:inclusion}) and $d\calS_{r+1}\Lambda^{k-1}\subset\calS_r\Lambda^k$ by the subcomplex property (\ref{eq:sr-subcpx}).
For the reverse containment, decompose the spaces as
\begin{align*}
d\calS_{r+1}\Lambda^{k-1} & = d\calP_{r+1}\Lambda^{k-1}+d\calJ_{r+1}\Lambda^{k-1} \\
\calS_r^-\Lambda^k & = \calP_{r-1}\Lambda^k + \calJ_{r-1} \Lambda^k +  d \calJ_{r}\Lambda^{k-1} +\kappa\calP_{r-1}\Lambda^{k+1} +\kappa d\calJ_r\Lambda^k
\end{align*}
Observe that $\calP_r\Lambda^k =d\calP_{r+1}\Lambda^{k-1} \oplus \kappa\calP_{r-1}\Lambda^{k+1}$ by (\ref{eq:hr-decomp}) and (\ref{eq:pr-sum-of-hr}), $\calJ_r\Lambda^k = \kappa d\calJ_r\Lambda^k
$ by (\ref{eq:kdJ-J}), and $d \calJ_{r+1}\Lambda^{k-1}$ appears as a summand for $d\calS_{r+1}\Lambda^{k-1}$.
Thus, by (\ref{eq:srk-decomp}), $\calS_r\Lambda^k\subset \calS_r^-\Lambda^k + d\calS_{r+1}\Lambda^{k-1}$.
\end{proof}

\begin{theorem}[Exactness]
\label{thm:exact}
Let $n,r \geq 1$. The sequence 
\[ 0 \rightarrow \bbR \rightarrow \calS_r^-\Lambda^0 \rightarrow \calS_r^-\Lambda^1 \rightarrow \cdots \rightarrow \calS_r^-\Lambda^{n-1} \rightarrow \calS_r^-\Lambda^n \rightarrow 0  \]
is exact.
\end{theorem}

\begin{proof}
By Lemma \ref{lem:sr-alt-decomp}, part (i), we can rewrite the beginning of the sequence as
\[ 0 \rightarrow \bbR \rightarrow \calS_r\Lambda^0 \rightarrow \calS_r^-\Lambda^1\rightarrow\cdots \]
The sequence is exact at $\calS_r\Lambda^0$ since the incoming and outgoing maps at $\calS_r\Lambda^0$ are the same as those in (\ref{eq:srk-seq}), which is exact. 
For $k \geq 1$, we will show that
\[ \calS^-_r\Lambda^{k-1} \rightarrow \calS^-_r\Lambda^{k} \rightarrow \calS^-_r\Lambda^{k+1} \]
is exact at $\calS_r^-\Lambda^k$ directly. 

Let $\omega \in \calS^-_r\Lambda^{k}$ and assume $d\omega = 0$. Using (\ref{eq:srm-directsumdecomp}), write $\omega = \sum_{i=1}^3 \omega_i$ where
\[ \omega_1 \in \calP_{r-1}^-\Lambda^k, \hspace{.2in} \omega_2 \in \calJ_r\Lambda^k, \hspace{.2in} \omega_3 \in d\calJ_r\Lambda^{k-1}. \]
Thus $d(\omega_3) = 0$ and $d(\omega_1 + \omega_2) = 0$. Since $\omega_1$ has maximum polynomial degree $r$ and $\omega_2$ has minimum polynomial degree $r+1$, we see that $d(\omega_1) = d(\omega_2) = 0$. 

Recall from (\ref{eq:prm-seq}) that $\calP_{r}^-\Lambda^{\bs}$ is exact. Thus, there exits $\mu_1 \in \calP_{r-1}^-\Lambda^{k-1} \subset  \calS^-_r\Lambda^{k-1}$ such that $d(\mu_1) = \omega_1$ (in particular, $\kappa(\omega_1)$ with an appropriate coefficient suffices).
Since $\omega_2 \in \calJ_r\Lambda^k$, we can write $\omega_2 = \kappa \mu_2$ for some polynomial $k+1$-form $\mu_2$. 
By hypothesis, $d(\kappa \mu_2) = d(\omega_2) =  0$, but $d$ is injective on the range of $\kappa$ by (\ref{eq:magic}). 
Therefore, $\kappa \mu_2 = \omega_2 = 0$. 
Also, since $\omega_3 \in d\calJ_r\Lambda^{k-1}$, we can write $\omega_3 = d\mu_3$ where $\mu_3 \in \calJ_r\Lambda^{k-1} \subset \calS^-_r\Lambda^{k-1}$, by (\ref{eq:srm-directsumdecomp}). 
Setting $\mu := \mu_1 + \mu_3 \in \calS^-_r\Lambda^{k-1}$, we have $d\mu = \omega$. 
\end{proof}

The $\calS_r^-\Lambda^k$ spaces also have a trace property analogous to the $\calS_r\Lambda^k$ spaces.
Recall that the \textit{trace} of a differential $k$-form on a codimension 1 hyperplane $f\subset\bbR^n$ is the pullback of the form via the inclusion map $f\hookrightarrow\bbR^n$.
Let $x^\alpha dx_\sigma$ be a form monomial as in (\ref{eq:form-mon}) and let $f$ be the hyperplane defined by $x_i=c$ for some fixed $1\leq i\leq n$ and  constant $c$.
Then
\[\tr_f (x^\alpha dx_\sigma)=\begin{cases} 0, & i\in\sigma,\\ \left(x^\alpha|_{x_i=c}\right)dx_\sigma, & i\not\in\sigma.
\end{cases}\]

\begin{theorem}[Trace property]
\label{thm:trace}
Let $n,r\geq 1$, $0\leq k\leq n$ and let $f$ be a hyperplane of $\bbR^n$ obtained by fixing one coordinate.  Then
\begin{equation}
\label{eq:trace}
\tr_f\calS_r^-\Lambda^k(\bbR^n) \subset \calS_{r}^-\Lambda^k(f).
\end{equation}
\end{theorem}
\begin{proof}
We use the result of~\cite[Theorem 3.5]{AA2014} and techniques from its proof to derive the result.
For a fixed constant $c \in \bbR$, set $f = \{ x \in \bbR^n~:~x_1 = c\}$. Using (\ref{eq:srm-directsumdecomp}), we need to show that the traces of $\calP_{r-1}^-\Lambda^k(\bbR^n)$, $\calJ_r\Lambda^k(\bbR^n)$, and $d\calJ_r\Lambda^{k-1}(\bbR^n)$ lie in
\[ \calS_r^-\Lambda^k(f) = \calS_{r-1}\Lambda^k(f) + \kappa \calS_{r-1}\Lambda^{k+1}(f) = \calP_{r-1}^-\Lambda^k(f) \oplus \calJ_r\Lambda^k(f) \oplus d\calJ_r\Lambda^{k-1}(f). \]
By \cite[Section 3.6]{AFW2006}, $\tr_f\calP_{r-1}^-\Lambda^k(\bbR^n) \subset \calP_{r-1}^-\Lambda^k(f)$ and by~\cite[Theorem 3.5]{AA2014}, 
\[ \tr_f d\calJ_r\Lambda^{k-1}(\bbR^n) = d\tr_f\calJ_r\Lambda^{k-1}(\bbR^n) \subset d\calS_r\Lambda^{k-1}(f) \subset \calS_{r-1}\Lambda^k(f). \]

It remains to show that $\tr_f\calJ_r\Lambda^k(\bbR^n) \subset \calS_r^-\Lambda^k(f)$.
Let $m$ be a $(k+1)$--form monomial $m$ with $\deg m\geq r$ and $\deg m - \ldeg~m \leq r - 1$.
By Proposition~\ref{prop:jrk-char}, it suffices to show that $\tr_f\kappa m \in \calS_r^-\Lambda^k(f)$.
Without loss of generality, we will write
\[ m  = x^\alpha dx_\sigma = x_1^{\alpha_1}x^{\alpha'}dx_{\sigma}, \]
where $\alpha':=\alpha-(\alpha_1,0,\cdots,0)$. 
Now, as the proof of~\cite[Theorem 3.5]{AA2014}, we break into cases according to whether or not $1\in\sigma$.
If $1\not\in\sigma$, we define $z$ by
\[\tr_f\kappa m= \tr_f\kappa \left(x_1^{\alpha_1}x^{\alpha'}dx_{\sigma}\right) = \kappa \left(c^{\alpha_1}x^{\alpha'}dx_{\sigma}\right) =:\kappa z.\]
If $\deg\kappa z\leq r$, then $\kappa z\in\calP_{r-1}^-\Lambda^k(f)\subset\calS_r^-\Lambda^k(f)$ and we are done.
So we presume $\deg\kappa z\geq r+1$, whence $\deg z\geq r$.  
If $\alpha_1\not=1$, $\deg z\leq \deg m$ and $\ldeg~z=\ldeg~m$.
If $\alpha_1=1$, $\deg z=\deg m -1$ and $\ldeg~z=\ldeg~m -1$.
Either way, $\deg z-\ldeg~z\leq \deg m -\ldeg~m\leq r-1$.
Thus, $\kappa z \in \calJ_r\Lambda^k(f)\subset \calS_r^-\Lambda^k(f)$, by the characterization of $\calJ_r\Lambda^k(f)$ from Proposition~~\ref{prop:jrk-char} and by the direct sum decomposition (\ref{eq:srm-directsumdecomp}), respectively.

Keeping $m$ as above, we now address the case $1\in\sigma$.
Define $w$ by
\[\tr_f\kappa m= \pm~c^{\alpha_1 + 1}x^{\alpha'}dx_{\tau}=:w\]
where $\tau \subset \{2, 3, \ldots, n\}$, and $\{1\} \cup \tau = \sigma$.
The sign of $w$ depends on the parity of the number of permutations required to reorder $\sigma(1),\cdots,\sigma(k+1)$ into $1,\tau(1),\cdots,\tau(k)$.
If $\deg w < r$, then $w \in \calP_{r-1}\Lambda^k(f) \subset \calS_r^-\Lambda^k(f)$.
So we presume $\deg w\geq r$.

We will show that $d \kappa w$, $\kappa d w \in \calS_r^-\Lambda^k(f)$, which by~(\ref{eq:magic}) implies that $w\in\calS_r^-\Lambda^k(f)$. 
Note that $\ldeg~w=\ldeg~m$, since $1\in\sigma$, and $\deg w\leq \deg m$ by definition of $\alpha'$.
Thus, $\deg w-\ldeg~w\leq \deg m-\ldeg~m\leq r-1$.
By Proposition~\ref{prop:jrk-char}, $\kappa w\in \calJ_r\Lambda^{k-1}(f)$ and thus $d\kappa w \in d\calJ_r\Lambda^{k-1}(f) \subset \calS_r^-\Lambda^k(f)$.

We split into cases one last time based on the inequality $\deg w\geq r$. 
If $\deg w = r$, then $\deg\kappa d w \leq r$ and we have $\kappa d w \in\calP_{r-1}^-\Lambda^k(f)\subset\calS_r^-\Lambda^k(f)$.
If $\deg w >r$, we have $\deg d w \geq r$ and $\deg dw=\deg w -1 \leq\deg m-1$.  
Since $d$ either preserves the linear degree of a form monomial or decreases it by one, we have $\ldeg~d w\geq \ldeg~w -1=\ldeg~m -1$.
Thus $\deg dw - \ldeg~d w \leq\deg m-\ldeg~m\leq r-1$.
Again by Proposition~\ref{prop:jrk-char}, $\kappa dw\in\calJ_r\Lambda^k(f)\subset \calS_r^-\Lambda^k(f)$.
\end{proof}

We now compute the dimension of $\calS_r^-\Lambda^k(\bbR_n)$ from the direct sum decomposition~(\ref{eq:srm-directsumdecomp}).
The computation of $\dim\calS_r\Lambda^k(\bbR^n)$ in~\cite{AA2014} does not rely on its direct sum decomposition~(\ref{eq:srk-decomp}) and, in particular, no formula for $\dim\calJ_r\Lambda^k(\bbR^n)$ is provided. 
We now derive such a formula.
We use $\dim X$ and $|X|$ interchangeably to denote the dimension of $X$ as a vector space over $\bbR$.

\begin{lemma}
\label{lem:dim-Jrk}
Fix $n\geq 1$.
For $r\geq 1$, $0\leq k\leq n$, we have 
\begin{equation}
\dim\calJ_r\Lambda^k(\bbR^n)  = \sum_{i=0}^k(-1)^i \left(A-B\right),
\label{eq:dim-Jrk}
\end{equation}
where
\begin{align*}
A := &~~ \sum_{d=k-i}^{\min\{n,\lfloor (r+i)/2\rfloor+k-i\}}2^{n-d}{n\choose d}{r-d+2k-i\choose d}{d\choose k-i}, \\[2mm]
B := &~~ \sum_{j=0}^{r+i}{n+j-1\choose j}{n\choose k-i}.
\end{align*}
\end{lemma}
\begin{proof}
Observe that
\[d\calS_r\Lambda^k = d\calP_r\Lambda^k \oplus d\calJ_r\Lambda^k = d\kappa\calP_{r-1}\Lambda^{k+1}\oplus d\calJ_r\Lambda^k.\]
Since $d$ is injective on the range of $\kappa$, we have $|d\kappa \calP_{r-1}\Lambda^{k+1}|= |\kappa \calP_{r-1}\Lambda^{k+1}|$ and $|d\calJ_r\Lambda^k|=|\calJ_r\Lambda^k|$.
Now, recall from (\ref{eq:srk-seq}) that $\calS_{r-\bs}\Lambda^\bs$ is exact.
Thus,
\begin{align}
|\calS_r\Lambda^k| &= |d\calS_r\Lambda^k| + |d\calS_{r+1}\Lambda^{k-1}| \label{eq:dim-srk-alt}\\
& = |\kappa\calP_{r-1}\Lambda^{k+1}| + |\calJ_r\Lambda^k| + |\kappa\calP_{r}\Lambda^{k}| + |\calJ_{r+1}\Lambda^{k-1}|. \notag
\end{align}
By (\ref{eq:pr-sum-of-hr}) and~\cite[Equation (3.14)]{AFW2006}, we have
\begin{align*}
|\kappa\calP_r\Lambda^k|+|\kappa\calP_{r-1}\Lambda^{k+1}| = \sum_{j=0}^r|\kappa\calH_j\Lambda^k|+|\kappa\calH_{j-1}\Lambda^{k+1}|
=\sum_{j=0}^r{n+j-1\choose j}{n\choose k}.
\end{align*}
Define $j_{r,k}:=|\calJ_r\Lambda^k(\bbR^n)|$ and $f_{r,k}:=j_{r,k}+j_{r+1,k-1}$ for ease of notation.
Using (\ref{eq:dim-srk-alt}) and the formula for $|\calS_r\Lambda^k(\bbR^n)|$ given in~\cite{AA2014}, we have
\begin{align}
f_{r,k} &= |\calS_r\Lambda^k| - \left(|\kappa\calP_{r}\Lambda^{k}| + |\kappa\calP_{r-1}\Lambda^{k+1}|  \right) \label{eq:frk-dim}\\
&= \left(\sum_{d=k}^{\min\{n,\lfloor r/2\rfloor+k\}}2^{n-d}{n\choose d}{r-d+2k\choose d}{d\choose k}\right) -\left(\sum_{j=0}^r{n+j-1\choose j}{n\choose k}\right).\notag
\end{align}
We can write $j_{r,k}$ as the telescoping sum
\begin{equation}
\label{eq:jrk-frk}
j_{r,k} = \sum_{i=0}^k(-1)^i f_{r+i,k-i}.
\end{equation}
Using (\ref{eq:jrk-frk}) with (\ref{eq:frk-dim}), we produce the formula in (\ref{eq:dim-Jrk}).
\end{proof}

\begin{table}[t]
\begin{centering}
\begin{tabular}{rr|rrrrrrr}
\text{} & \\[-2mm]
& $k$ & $r$=1 & 2 & 3 & 4 & 5 & 6 & 7 \\
\text{} & \\[-2mm]
\hline
\text{} & \\[-2mm]
$n$=1 & 0 & 2 & 3 & 4 & 5 & 6 & 7 & 8 \\
    & 1 & 1 & 2 & 3 & 4 & 5 & 6 & 7 \\[2mm]
\hline
\text{} & \\[-2mm]
$n$=2 & 0 & 4 & 8  & 12 & 17 & 23 & 30 & 38 \\ 
    & 1 & 4 & 10 & 17 & 26 & 37 & 50 & 65 \\ 
    & 2 & 1 & 3  & 6  & 10 & 15 & 21 & 28 \\[2mm] 
\hline 
\text{} & \\[-2mm]
$n$=3 & 0 & 8  & 20 & 32 & 50  & 74  & 105 & 144 \\ 
    & 1 & 12 & 36 & 66 & 111 & 173 & 255 & 360 \\ 
    & 2 &  6 & 21 & 45 & 82  & 135 & 207 & 301 \\ 
    & 3 &  1 & 4  & 10 & 20  & 35  & 56  & 84  \\[2mm] 
\hline 
\text{} & \\[-2mm]   
$n$=4 & 0 & 16 & 48  & 80  & 136 & 216 & 328  & 480  \\
    & 1 & 32 & 112 & 216 & 392 & 656 & 1036 & 1563 \\
    & 2 & 24 & 96  & 216 & 422 & 746 & 1227 & 1910 \\
    & 3 & 8  & 36  & 94  & 200 & 375 & 644  & 1036 \\
    & 4 & 1  & 5   & 15  & 35  & 70  & 126  & 210  \\[2mm] 
\hline
\end{tabular}\\[2mm]
\end{centering}
\caption{Dimension of $\calS_r^-\Lambda^k(\square_n)$ for $1\leq n\leq 4$, $0\leq k\leq n$, and $1\leq r\leq 7$, computed using Theorem~\ref{thm:dim-srmk}.}
\label{tab:srm-dim}
\end{table}

\begin{theorem}
\label{thm:dim-srmk}
Fix $n,r\geq 1$ and $0\leq k\leq n$.
Then
\begin{equation}
\label{eq:dim-Srmk}
\dim\calS_r^-\Lambda^k(\bbR^n) = \dim\calP_r^-\Lambda^k(\bbR^n) + \dim\calJ_r\Lambda^k(\bbR^n)  + \dim\calJ_{r} \Lambda^{k-1}(\bbR^n).
\end{equation}
Further, each summand in (\ref{eq:dim-Srmk}) has a closed-form expression in terms of binomial coefficients depending only on $n$, $k$, and $r$.
\end{theorem}
\begin{proof}
Again, since $d$ is injective on the range of $\kappa$, we have $|d\calJ_r\Lambda^{k-1}|=|\calJ_r\Lambda^{k-1}|$.
Using this with (\ref{eq:srm-directsumdecomp}), we can write
\[|\calS_r^-\Lambda^k| = |\calP_r^-\Lambda^k| + |\calJ_r\Lambda^k|  + |d \calJ_{r} \Lambda^{k-1}| = |\calP_r^-\Lambda^k| + |\calJ_r\Lambda^k|  + |\calJ_{r} \Lambda^{k-1}|.\]
From~\cite{AFW2006,AFW2010}, we have
\begin{equation}
\label{eq:dim-prm}
|\calP_r^-\Lambda^k| = {r+n\choose r+k}{r+k-1\choose k}.
\end{equation}
We have the requisite expressions for $|\calJ_r\Lambda^k|$ and $|\calJ_{r} \Lambda^{k-1}|$ from Lemma~\ref{lem:dim-Jrk}.
\end{proof}
We use Theorem~\ref{thm:dim-srmk} and Lemma~\ref{lem:dim-Jrk} to compute the dimension of $\calS_r^-\Lambda^k(\square_n)$ for $1\leq n\leq 4$, $0\leq k\leq n$, and $1\leq r\leq 7$ and report the results in Table~\ref{tab:srm-dim}.

\section{Degrees of Freedom, Unisolvence, and Minimality}
\label{sec:unisolv-min}

We now state and count a set of degrees of freedom associated to $\calS_r^-\Lambda^k(\square_n)$.
The degrees of freedom associated to a $d$-dimensional sub-face $f$ of $\square_n$ are
\begin{equation}
\label{eq:Srm-dofs}
u\longmapsto \int_{f}(\tr_f\,u)\wedge q, \quad q\in \calP_{\,r-2(d-k)-1}\Lambda^{d-k}(f)~\oplus~d\calH_{r-2(d-k)+1}\Lambda^{d-k-1}(f),
\end{equation}
for any $k\leq d\leq\min\{n,\lfloor r/2\rfloor+k\}$.
Observe that the first summand of the indexing space is the indexing space for $\calS_{r-1}\Lambda^k(f)$, reflecting the fact that $\calS_r^-\Lambda^k\supset\calS_{r-1}\Lambda^k$. 
The sum is direct since $d\calH_{r-2(d-k)+1}\Lambda^{d-k-1}\subset\calH_{r-2(d-k)}\Lambda^{d-k}$.
The dimension of $\calP_r\Lambda^k(\bbR^n)$ is given (see e.g.~\cite{AFW2006}) by
\begin{equation}
\label{eq:prk-dim}
\dim\calP_r\Lambda^k(\bbR^n) = {r+n\choose r+k}{r+k\choose k}.
\end{equation}
Applying (\ref{eq:prk-dim}), we have that
\[\dim\calP_{\,r-2(d-k)-1}\Lambda^{d-k}(f)={r-d+2k-1 \choose r-d+k-1}{r-d+k-1\choose d-k}.\]
It is shown in~\cite[Theorem 3.3]{AFW2006} that
\[\dim d\calH_{r+1}\Lambda^{k-1}(\bbR^n)  =\dim\kappa\calH_r \Lambda^k(\bbR^n) = {n+r\choose n-k}{r+k-1\choose k-1},\]
and thus
\[\dim d\calH_{r-2d+2k+1}\Lambda^{d-k-1}(f) = {r-d+2k\choose k}{r-d+k-1\choose d-k-1}. \]
Note that when $k=d$, we have $d\calH_{r+1}\Lambda^{-1}(f)=d(0)=0$, so the dimension is zero.  The above formula remains valid if we interpret $\ds{r-1\choose -1}$ as 0.
There are $2^{n-d}{n\choose d}$ $d$-dimensional faces of $\square_n$ so the total number of degrees of freedom in (\ref{eq:Srm-dofs}) is
\begin{align}
\sum_{d=k}^{\min\{n,\lfloor r/2\rfloor+k\}}2^{n-d}{n\choose d}\left({r-d+2k-1 \choose r-d+k-1}\right.  & {r-d+k-1\choose d-k}\notag \\
 & + \left. {r-d+2k\choose k}{r-d+k-1\choose d-k-1}\right).
\label{eq:Srm-dofs-count}
\end{align}
To prove that the degrees of freedom in (\ref{eq:Srm-dofs}) are unisolvent for $\calS_r^-\Lambda^k(\square_n)$ we will need to consider the subspace of $\calS_r^-\Lambda^k(\square_n)$ that has vanishing trace on $\p\square_n$.
For this, we will use the notation
\[\calS_r^-\Lambda^k_0(\square_n):=\left\{~\omega\in\calS_r^-\Lambda^k(\square_n)~:~\tr_f\omega=0~\text{for every $(n-1)$-subface of $\square_n$}~\right\}.\]
The next result is the analogue of~\cite[Proposition 3.7]{AA2014} for the $\calS_r^-\Lambda^k(\square_n)$ family.
\begin{lemma}
\label{lem:vanish-trace}
If $\omega\in\calS_r^-\Lambda^k_0(\square_n)$ and
\begin{align}
\int_{\square_n} \omega\wedge p = 0, &\qquad p\in P_{\,r-2(n-k)-1}\Lambda^{n-k}(\square_n) \label{eq:vanish-p}\\[2mm]
\int_{\square_n} \omega\wedge dh = 0, &\qquad h\in 
\calH_{r-2(n-k)+1}\Lambda^{n-k-1}(\square_n) \label{eq:vanish-h}
\end{align}
then $\omega\equiv 0$.
\end{lemma}

\begin{proof}
Let $\omega\in\calS_r^-\Lambda^k_0$.
By the subcomplex property (\ref{eq:subcpx}), we have that $d\omega\in\calS_{r}^-\Lambda^{k+1}$.
Recalling the definition  $\calS_r^-\Lambda^{k+1}=\calS_{r-1}\Lambda^{k+1} +\kappa\calS_{r-1}\Lambda^{k+2}$ and the fact that $d$ is injective on the range of $\kappa$, we have that $d\omega\in\calS_{r-1}\Lambda^{k+1}$, a non-trimmed serendipity space.
Let $f$ be any $(n-1)$--face of $\square_n$ and recall that $d$ commutes with $\tr_f$.  
Thus, $\tr_f d\omega = d\tr_f\omega = 0$, meaning $d\omega\in\calS_{r-1}\Lambda^{k+1}_0$.

By Stokes' theorem, we have
\[\int_{\square_n} d\omega\wedge\mu = \pm\int_{\square_n} \omega\wedge d\mu,\qquad \mu\in\Lambda^{n-k-1}(\square_n).\]
Suppose $\mu\in\calP_{r-2(n-k)+1}\Lambda^{n-k-1}$ so that
\[d\mu\in \calP_{\,r-2(n-k)-1}\Lambda^{n-k}~\oplus~d\calH_{r-2(n-k)+1}\Lambda^{n-k-1}.\]
By (\ref{eq:vanish-p}) and (\ref{eq:vanish-h}),  $\int \omega\wedge d\mu$ vanishes for all such $\mu$ and by the above equation $\int d\omega\wedge\mu$ vanishes for all such $\mu$ as well.
Thus, by~\cite[Proposition 3.7]{AA2014} with $r$ and $k$ replaced by $r-1$ and $k+1$, respectively, we have $d\omega=0$.

By Theorem~\ref{thm:srm-directsumdecomp}, we can write $\omega=d\alpha+\kappa \beta$ where $d\alpha\in\calS_{r-1}\Lambda^k$ and $\kappa\beta\in\calS_r\Lambda^k$.
Since $d\omega=0$ and $d$ is injective on the range of $\kappa$, we must have $\kappa\beta=0$.
Thus $\omega=d\alpha\in\calS_{r-1}\Lambda^k$.
Since (\ref{eq:vanish-p}) holds, we can apply~\cite[Proposition 3.7]{AA2014} with $r$ replaced by $r-1$ to conclude that $\omega\equiv0$. 
\end{proof}

We can now establish unisolvence in the classical sense, namely, that an element $u\in\calS_r^-\Lambda^k$ is uniquely determined by the values of the degrees of freedom applied to $u$.

\newpage

\begin{theorem}[Unisolvence]
\label{thm:unisolv}
If $u\in\calS_r^-\Lambda^k(\square_n)$ and all the degrees of freedom in (\ref{eq:Srm-dofs}) vanish then $u\equiv 0$.
\end{theorem}
%
%

\begin{proof}
We use induction on $n$. 
The base case $n=1$ is trivial.
Let $\omega\in\calS_r^-\Lambda^k(\square_n)$ such that all the degrees of freedom in (\ref{eq:Srm-dofs}) vanish.
On a face $f$ of dimension $n-1$, $\tr_f\omega\in\calS_r^-\Lambda^k(f)$ by the trace property~(\ref{eq:trace}).
Since all the degrees of freedom for $\tr_f\omega$ vanish, $\tr_f\omega\equiv 0$ by the inductive hypothesis.
Thus, $\omega\in\calS_r^-\Lambda^k_0(\square_n)$.
By Lemma~\ref{lem:vanish-trace}, $\omega\equiv 0$.
\end{proof}

The careful combinatorial argument carried out in~\cite{AA2014} to establish unisolvence for the $\calS_r\Lambda^k$ spaces and their associated degrees of freedom is essential to the proof of unisolvence for the $\calS_r^-\Lambda^k$ spaces just given, as it is invoked at the end of the proof of Lemma~\ref{lem:vanish-trace}.
Notably, our proof did not require the dimension of $\calS_r^-\Lambda^k(\square_n)$ to equal the associated number of degrees of freedom as a hypothesis.
We examine this point further in the discussion of future directions at the end of the paper and in Appendix~\ref{appdx:dim-eq-exs}.

We now turn to the topic of the minimality of the $\calS_r^-\Lambda^k$ spaces.
For this, we will employ the theory of \textit{finite element systems}, developed and applied by Christiansen and collaborators in~\cite{Chr08M3AS,Chr09AWM,CG2015,ChrMunOwr11,CR2016}.
We will not redefine the full framework here as we are interested only in a very specific context, similar to the examples studied in~\cite{CG2015}.

We have shown that the $\calS_r^-\Lambda^k$ spaces have the subcomplex and trace properties in Theorems \ref{thm:subcpx} and \ref{thm:trace}, respectively.
These properties ensure that the collection of spaces $\{\calS_r^-\Lambda^0(\square_n),\ldots,\calS_r^-\Lambda^n(\square_n)\}$ constitute a finite element system, for any fixed $n,r\geq 1$.
Since the associated augmented co-chain complex for this sequence was shown to be exact in Theorem~\ref{thm:exact}, the system is said to be \textit{locally exact}.
Whenever unisolvence holds in the sense established in Theorem~\ref{thm:unisolv} \textit{and} the number of degrees of freedom equals the dimension of the associated trimmed serendipity spaces in the sequence, the system is said to \textit{admit extensions} and be \textit{compatible}.
In such cases, we can apply the following result, specialized to the case of cubical meshes.

\begin{lemma}[{\cite[Corollary 3.2]{CG2015}}]
\label{lem:minimality}
Suppose that $A$ is a finite element system on $\square_n$ and that $B$ is a compatible finite element system containing $A$.
Suppose that
\begin{equation}
\label{eq:mcfes}
\dim B^k_0(\square_n) = \dim A^k_0(\square_n) + \dim \rmH^{k+1}\left(A^\bs_0(\square_n)\right).\end{equation}
Then $B$ is minimal among compatible finite element systems containing $A$.
\end{lemma}

In (\ref{eq:mcfes}), $\rmH^{k+1}\left(A^\bs_0(\square_n)\right)$ denotes the $k+1$ homology group of the system $A^\bs_0$; the subscript $0$ again indicates vanishing trace on all $n-1$ dimensional subfaces.
Note that the system $A^\bs_0$ need not be locally exact and hence need not have vanishing homology. 
We can compute the dimension of the homology group by
\[\dim\left(\rmH^{k+1}\left(A^\bs_0(\square_n)\right)\right) =\dim \left(\ker d:A^{k+1}_0\raw A^{k+2}_0\right)-\dim (dA^{k}_0).\]
We apply the lemma as follows.

\begin{theorem}[Minimality]
\label{thm:minimal}
For $n=2$ and $n=3$, the system $\calS_r^-\Lambda^\bs(\square_n)$ is a minimal compatible finite element system containing $\calP_{r-1}\Lambda^\bs(\square_n)$.
\end{theorem}

\begin{remark}
Theorem~\ref{thm:minimal} is stated as applying only to dimensions $n=2$ and $n=3$, however, it holds in any setting for which the number of degrees of freedom equals the dimension of the associated trimmed serendipity spaces.
This includes at least all $r$ values up to 100 for $n=4$ and $n=5$.
We discuss this point further in Section~\ref{sec:conc} and Appendix~\ref{appdx:dim-eq-exs}.
\end{remark}

\begin{proof}
We set $A^k(\square_n):=\calP_{r-1}\Lambda^k(\square_n)$ and $B^k=\calS_r^-\Lambda^k(\square_n)$ and show that Lemma~\ref{lem:minimality} applies.
Note that $\calP_{r-1}\Lambda^\bs(\square_n)$ is a non-compatible finite element system as it satisfies the subcomplex and trace properties but is not locally exact.
We have already discussed why $\calS_r^-\Lambda^k(\square_n)$ is a compatible finite element system and shown that ${\calP_{r-1}\Lambda^k(\square_n) \subset \calS_r^-\Lambda^k(\square_n)}$.
By~\cite[Proposition 4.5]{CG2015}, 
\[\dim\calP_{r-1}\Lambda^k_0(\square_n)=\dim\calP_{r-2(n-k)-1}\Lambda^{n-k}(\square_n).\]
In the proof of~\cite[Lemma 4.13]{CG2015}, it is shown that
\[\dim \rmH^{k}\left(\calP_r\Lambda^\bs_0(\square_n)\right)=\dim\kappa\calH_{r+2k-2n-1}\Lambda^{n-k+1}(\square_n).\]
By~\cite[Theorem 3.3]{AFW2006},
\[\dim\kappa\calH_{r+2k-2n-1} \Lambda^{n-k+1}(\square_n)=\dim d\calH_{r+2k-2n}\Lambda^{n-k}(\square_n).\]
Replacing $r$ by $r-1$ and $k$ by $k+1$, we have 
\[\dim \rmH^{k+1}\left(\calP_{r-1}\Lambda^\bs_0(\square_n)\right)=\dim d\calH_{r+2k-2n+1}\Lambda^{n-k-1}(\square_n).\]
Applying Theorem~\ref{thm:unisolv}, we have
\begin{align*}
\dim \calS_r^-\Lambda^k_0(\square_n) & = \text{\# of degrees of freedom associated to the interior of $\square_n$} \\
&= \dim \calP_{\,r-2(n-k)-1}\Lambda^{n-k}(\square_n)~+~d\calH_{r-2(n-k)+1}\Lambda^{n-k-1}(\square_n) \\
&= \dim\calP_{r-1}\Lambda^k_0(\square_n) + \dim \rmH^{k}\left(\calP_r\Lambda^\bs_0(\square_n)\right).
\end{align*}
Therefore, Lemma~\ref{lem:minimality} applies and minimality is proved.
\end{proof}

We close this section with an examination of the computational benefit that using a minimal compatible finite element system can provide by comparing the use of trimmed serendipity elements in place of regular serendipity or tensor product (\nedelec) elements for a simple problem.
Consider the standard mixed formulation of the Dirichlet problem for the Poisson equation on a cubical domain $\Omega\subset\bbR^3$: Given $f$, find $\bfu\in H(\Div,\Omega)$ and $p\in L^2(\Omega)$ such that:
\begin{alignat}{2}
\int_\Omega \bfu\cdot\bfv &= \int_\Omega\Div\bfv~p, && \quad\quad\bfv\in H(\Div,\Omega), \label{eq:pois-1} \\
\int_\Omega \Div\bfu~q &= \int_\Omega f~q, &&\quad\quad q\in L^2(\Omega). \label{eq:pois-2}
\end{alignat}
We remark that (\ref{eq:pois-1})-(\ref{eq:pois-2}) is one instance of the Hodge Laplacian problem studied in finite element exterior calculus~\cite{AFW2006,AFW2010} and its analysis serves as the foundation for many applications, such as the movement of a fluid through porous media via Darcy flow~\cite{APWY2007}, diffusion via the heat equation~\cite{AC2017}, wave propagation~\cite{GHZ2017}, and various nonlinear partial differential equations.

A finite element method for (\ref{eq:pois-1})-(\ref{eq:pois-2}) is determined by selecting finite-dimensional subspaces $\Lambda^2_h\subset H(\Div,\Omega)$ and $\Lambda^3_h\subset L^2(\Omega)$ and solving the problem: find $\bfu_h\in \Lambda^2_h$ and $p_h\in \Lambda^3_h$ such that:
\begin{alignat}{2}
\int_\Omega \bfu_h\cdot\bfv_h &= \int_\Omega\Div\bfv_h~p_h, && \quad\quad\bfv_h\in \Lambda^2_h, \label{eq:pois-h-1} \\
\int_\Omega \Div\bfu_h~q_h &= \int_\Omega f~q_h, &&\quad\quad q_h\in \Lambda^3_h. \label{eq:pois-h-2}
\end{alignat}
Supposing that $\Omega$ is meshed by cubes, we compare three choices for the pair $(\Lambda^2_h,\Lambda^3_h)$ with at least $O(h^r)$ decay in the approximation of $p$, $\bfu$, and $\Div \bfu$ in the appropriate norms: tensor product elements $(\calQ_r^-\Lambda^2, \calQ_r^-\Lambda^3)$, serendipity elements $(\calS_r\Lambda^2, \calS_{r-1}\Lambda^3)$, and trimmed serendipity elements $(\calS_r^-\Lambda^2, \calS_r^-\Lambda^3)$.
We report the number of degrees of freedom associated to a single mesh element in Table~\ref{tab:mixed-dim-comp}.
As is evident from the table, the trimmed serendipity elements require the fewest degrees of freedom of any of the three choices.
Notably, the trimmed serendipity choice uses the same number of degrees of freedom as the tensor product elements in the lowest order case ($r=1$) while using strictly fewer than either of the other choices in all other cases.

\begin{table}
\begin{centering}
\begin{tabular}{r|rl|rl|rl|}
$r$ & \multicolumn{2}{c}{$|\calQ_r^-\Lambda^2|+|\calQ_r^-\Lambda^3|$} & \multicolumn{2}{c}{$|\calS_{r}\Lambda^2|+|\calS_{r-1}\Lambda^3|$} & \multicolumn{2}{c}{$|\calS_r^-\Lambda^2|+|\calS_r^-\Lambda^3|$} \\[1mm]
\hline
&&&&& \text{}\\[-2mm]
1 & 6+1 &= ~~~7     & 18+1 &= ~~19   & 6+1 &= ~~7 \\[1mm]
2 & 36+8 &= ~~44    & 39+4 &= ~~43   & 21+4 &= ~~25 \\[1mm]
3 & ~~108+27 &= ~~135\qquad 
                    & ~~72+10  &= ~~82\qquad 
                                     & 45+10 &= ~~55\qquad \\[1mm]
4 & 240+64 &= ~~304 & 120+20 &=~~140 & 82+20 &= ~~102\\[2mm]
\hline 
\end{tabular}\\[2mm]
\end{centering}
\caption{Comparison of dimension counts for a single cube element among pairs of tensor product, serendipity and trimmed serendipity spaces suitable for a mixed finite element formulation of the Poisson problem.}
\label{tab:mixed-dim-comp}
\end{table}

\section{Summary, Outlook, and Future Directions}
\label{sec:conc}

In this paper, we have defined spaces of trimmed serendipity finite element differential forms on $n$-dimensional cubes and demonstrated how their relation to the non-trimmed serendipity spaces are, in all essential ways, analogous to the relation of the trimmed and non-trimmed polynomial differential form spaces on simplices.
Accordingly, it is natural to treat them as a ``fifth column'' of the Periodic Table of Finite Elements~\cite{AL2014}.
The ease with which the trimmed serendipity spaces arise in the exterior calculus setting echoes the fact that instances of their vector calculus analogues have been discovered from the related but distinct frameworks of Arbogast and Correa~\cite{AC2015} and Cockburn and Fu~\cite{CF2016}, as detailed in Appendix~\ref{appdx:vec-calc}.

A minor point mentioned after the proof of Theorem~\ref{thm:unisolv} hints at an important direction for future research.
While we have shown that the degrees of freedom given in (\ref{eq:Srm-dofs}) are unisolvent for $\calS_r^-\Lambda^k(\square_n)$, this only establishes that the number of degrees of freedom is greater than or equal to the dimension of $\calS_r^-\Lambda^k(\square_n)$.
Using Mathematica, we verified that that formula (\ref{eq:Srm-dofs-count}) and the closed-form expression for $\dim\calS_r^-\Lambda^k(\square_n)$ from Theorem~\ref{thm:dim-srmk} are in fact equal for $1\leq n\leq 5$, $1\leq r \leq 100$, and $0\leq k\leq n$, covering many more than the cases of practical relevance to modern applications.
In the cases $n=2$ and $n=3$, we also confirmed by direct proof that the spaces have the equal dimension for any $r$; these proofs appear in Appendix~\ref{appdx:dim-eq-exs}.

A more promising approach toward the same goal is to construct a basis for $\calS_r^-\Lambda^k_0(\square_n)$, count its dimension, and sum over sub-faces of $\square_n$.
Such an approach is used by Arbogast and Correa~\cite{AC2015} in their study of mixed methods on quadrilaterals, but extending it to hexahedra, or even just to cubes, introduces significant additional subtleties regarding the linear independence of spanning sets of specific sets of polynomial differential forms.
We plan to explore this approach in future work, not only to establish this particular equality, but also to provide a practical computational basis so that the trimmed serendipity spaces can begin to see their benefits realized in practical application settings. 

\newpage

\text{}\\[0mm]
\textbf{Acknowledgements.}
AG was supported in part by NSF Award 1522289.  
The authors would like to thank the anonymous referees of the paper for their careful reading and helpful suggestions in their reviews.
The accepted version of this manuscript will appear in Mathematics of Computation, published by the American Mathematical Society.
\vspace{-3mm}

\bibliography{arxiv-trimmedSrdp}{}
\bibliographystyle{abbrv}

\appendix
\section{Trimmed serendipity spaces in vector calculus notation}
\label{appdx:vec-calc}

To characterize the relationship between the trimmed serendipity spaces of differential forms and finite element families described in traditional vector calculus notation, we will need some additional notation.
First we recall classical notation for spaces of polynomials and polynomial vector fields as used in~\cite{AC2015}.
Let $\mathbb{P}_r$ denote the space polynomials of degree at most $r$ and $\tilde{\mathbb{P}}_r$ the space of homogeneous polynomials of degree exactly $r$.
The number of variables (typically two or three) is implied from context.
For $n=2$, define
\[\mathbf{x}\tilde{\mathbb{P}}_r := \text{span}\left\{\twovec {x_1 p}{x_2p}~:~ p\in\tilde{\mathbb{P}}_r\right\}.\]
The above definition extends to $n=3$ by using $x_3p$ as the third component of the vector.
For $n=2$, define a ``2D curl'' operator on a scalar field $w$ as the gradient operator followed by a rotation of $\pi/2$ clockwise, i.e.
\[\curl\, w := \rot\,\nabla w = \begin{bmatrix} 0 & {1} \\ {-1} & 0\end{bmatrix}\twovec{\p w/\p x_1}{\p w/\p x_2} = \twovec{\p w/\p x_2}{-\p w/\p x_1}.\]
Thus, we recover the statement $\Div\curl\,w=0$ for any $w\in C^2$.

To convert a vector field to its corresponding differential form, we use the \textit{flat} operator, $\flat$, following the conventions of Abraham et al.~\cite{AMR1988}; see also Hirani~\cite{H2003}.
Given a scalar field $w$ on $\bbR^2$, there is an associated 0-form $w^\flat:=w$ and an associated 2-form $w^\flat:=w\,dx_1dx_2$.
It will be clear from context whether $w^\flat$ should be interpreted as a 0-form or a 2-form.
Given a vector field $\bfv=\twovec{v_1}{v_2}$ on $\bbR^2$, define $\bfv^\flat:=v_1dx_1+v_2dx_2$ and $\rot\left(\bfv^\flat\right):=\left(\rot\bfv\right)^\flat$.
In $\bbR^3$, given a scalar field $w$ on $\bbR^3$, there is an associated 0-form $w^\flat:=w$ and an associated 3-form $w^\flat:=w\,dx_1dx_2dx_3$.
Given a vector field $\bfv=\threevec{v_1}{v_2}{v_3}$ on $\bbR^3$, the associated 1-form is defined by $\bfv^\flat:=v_1dx_1+v_2dx_2+v_3dx_3$ and the associated 2-form by $\bfv^\flat:=v_1dx_2dx_3-v_2dx_1dx_3+v_3dx_1dx_2$.
Again, the kind of flat operator to be used will be obvious from context.
\\

\paragraph{Proof of Proposition~\ref{prop:ac-comp}.}

The space $AC_r(\hat E)$ from~\cite{AC2015} refers to a pair of spaces $\left(\mathbf{V}_{AC}^r,W_{AC}^r\right)\subset H(\div,\hat E)\times L^2(\hat E)$, where $\hat E=\square_2=[-1,1]^2$ is the reference element.
These spaces are defined to be
\[\mathbf{V}_{AC}^r := {\mathbb{P}}_r^2\oplus \mathbf{x}\tilde{\mathbb{P}}_r \oplus\mathbb{S}_r,\qquad\text{and}\qquad W_{AC}^r:=\mathbb{P}_r.\]
Since $L^2(\hat E)$ here corresponds to $\Lambda^2(\square_2)$ in finite element exterior calculus notation, we observe that $\left(W_{AC}^r\right)^\flat:=\left(\mathbb{P}_r\right)^\flat=\calP_r\Lambda^2(\square_2)$.
Further, we have $\calP_r\Lambda^2(\square_2)=\calS_{r}\Lambda^2(\square_2)$, which is identical to $\calS_{r+1}^-\Lambda^2(\square_2)$ by Lemma~\ref{lem:sr-alt-decomp}.

Turning to $\mathbf{V}_{AC}^r$, observe that $\left({\mathbb{P}}_r^2\right)^\flat=\calP_r\Lambda^1(\square_2)$.
Note that $\rot$ is an automorphism on $\calP_r\Lambda^1(\square_2)$, i.e.\ $\rot\calP_r\Lambda^1(\square_2)=\calP_r\Lambda^1(\square_2)$.
Using $x$ and $y$ in place of $x_1$ and $x_2$ and omitting `span' notation for ease of reading, we also see that
\begin{align*}
\rot\left(\mathbf{x}\tilde{\mathbb{P}}_r\right)^\flat & = \left(\left\{\rot\twovec {xp}{y p}~:~ p\in\tilde{\mathbb{P}}_r\right\}\right)^\flat = \left(\left\{\twovec {yp}{-x p}~:~ p\in\tilde{\mathbb{P}}_r\right\}\right)^\flat \\
 &= \left\{-\kappa p\,dxdy~:~p\in\calH_r\Lambda^0(\square_2)\right\} = \kappa\calH_r\Lambda^0(\square_2).
\end{align*}
Hence, we have $\rot\left({\mathbb{P}}_r^2\oplus\mathbf{x}\tilde{\mathbb{P}}_r\right)^\flat=\calP_r\Lambda^1(\square_2)\oplus\kappa\calH_r\Lambda^0(\square_2)=\calP_{r+1}^-\Lambda^1(\square_2)$.

The space $\mathbb{S}_r$ is a space of ``supplemental'' vectors that satisfy certain conditions described in~\cite{AC2015}.
Foremost, the space $\mathbb{S}_r$ satisfies the containment
\[\mathbb{S}_r\subset\curl\left(\mathbb{P}_{r+1}\oplus\text{span}\left\{x^iy^j~:~i+j=r+2,~1\leq i,j\leq r+1\right\}\right),\]
Further, the elements of $\mathbb{S}_r$ are required to have normal components on $\square_2$ that are polynomials of degree $r$.
The authors present the following basis for $\mathbb{S}_r$ for $r\geq 1$:
\begin{align*}
\text{basis for $\mathbb{S}_r$} 
 &= \left\{\curl(x^{r-1}(1-x^2)y),~~\curl(xy^{r-1}(1-y^2))\right\} \\
 &= \left\{\twovec {x^{r-1}(1-x)^2 }{ -y\left((r-1)x^{r-2} -(r+1) x^r\right)}, \twovec {x \left((r-1)y^{r-2} -(r+1) y^r\right) }{ -y^{r-1}(1-y^2)} \right\} \\
 &=: \{\hat\sigma_1,\hat\sigma_2\}
\end{align*}
Looking at the homogeneous degree $r+1$ part of $\hat\sigma_1$, we see that
\[{\hat\sigma_1}^\flat= \twovec{x^{r+1}}{(r+1)x^ry}^\flat+\bfv^\flat=x^{r+1}\, dx+(r+1)x^ry\,dy +\bfv^\flat.\]
Applying $\rot$ to both sides, we have that
\[\rot{\hat\sigma_1}^\flat= (r+1)x^ry\,dx - x^{r+1}\, dy +\rot\bfv^\flat = d\kappa (x^ry\,dx) + \rot\bfv^\flat. \]
Note that $d\kappa (x^ry\,dx) \in d\kappa\calH_{r+1,1}\Lambda^1(\square_2)=d\calJ_{r+1}\Lambda^{0}(\square_2)$\footnote{Recall that $\calH_{r,\ell}\Lambda^k(\bbR^n)=0$ if $\ell>\min(r,n-k)$; the relevant case here is $\ell>\min(r,2-1)=1$.} and $\rot\bfv^\flat\in\calP_r\Lambda^1(\square_2)$.
Hence, $\rot\hat\sigma_1^\flat\in d\calJ_{r+1}\Lambda^{0}(\square_2)+\calP_r\Lambda^1(\square_2)\subset\calS_{r+1}^-\Lambda^1(\square_2)$ and similarly, $\rot\hat\sigma_2^\flat\in\calS_{r+1}^-\Lambda^1(\square_2)$.
Observe that $\rot \hat\sigma_1^\flat$ and $\rot\hat\sigma_2^\flat$ are linearly independent and have distinct, non-zero projections on to $d\calJ_{r+1}\Lambda^{0}(\square_2)$.
Thus, given a basis $\{\bfv_1,\ldots,\bfv_m\}$ for ${\mathbb{P}}_r^2\oplus\mathbf{x}\tilde{\mathbb{P}}_r$, the set $\{\bfv_1,\ldots,\bfv_m,\hat\sigma_1,\hat\sigma_2\}$ is a basis for $\mathbf{V}_{AC}^r$ and the set
\[\{\rot\bfv_1^\flat,\ldots,\rot\bfv_m^\flat,\rot\hat\sigma_1^\flat,\rot\hat\sigma_2^\flat\}\]
is a basis for $\calS_{r+1}^-\Lambda^1(\square_2)$.
This proves Proposition~\ref{prop:ac-comp}.\\

\paragraph{Proof of Proposition~\ref{prop:cf-comp}.}
We now turn to the paper by Cockburn and Fu~\cite{CF2016}.
Rather than restate all their definitions, we translate their notation to the Arbogast-Correa notation or finite element exterior calculus notation as we analyze their spaces.
First, we look at the sequence $S_{2,r}^\square(K)$ from their Theorem 3.3.
Applying the flat operator to the first space, we get
\[
\left(\mathscr P_{r+1}(x,y)\oplus\delta H^{2,I}_{r+1}\right)^\flat = \calP_{r+1}\Lambda^0(\square_2)\oplus \left(\text{span}\left\{xy^{r+1}, x^{r+1}y\right\}\right)^\flat=\calS_{r+1}\Lambda^0(\square_2).\]
Recall, by Lemma~\ref{lem:sr-alt-decomp}, that $\calS_{r+1}\Lambda^0(\square_2)=\calS_{r+1}^-\Lambda^0(\square_2)$.
The second space is written as a direct sum of three components.
Applying the flat operator to each, we find that
\begin{align*}
\left(\pmb{\mathscr{P}}_{r}(x,y)\right)^\flat & = \calP_{r}\Lambda^1(\square_2), \\
\left(\bfx \times \tilde{\mathscr{P}}_r(x,y)\right)^\flat &= \left(\rot\bfx \tilde{\mathbb{P}}_r\right)^\flat = \kappa\calH_r\Lambda^0(\square_2),\\
\left(\nabla\delta H^{2,I}_{r+1}\right)^\flat &= \left(\text{span}\nabla\left\{xy^{r+1},\, x^{r+1}y\right\}\right)^\flat = \text{span}\left\{d\kappa(xy^r dy),\, d\kappa(x^ry dx)\right\}\\
 & = d\kappa\calH_{r+1,1}\Lambda^1(\square_2) = d\calJ_{r+1}\Lambda^0(\square_2).
\end{align*}
By the direct sum decomposition (\ref{eq:srm-directsumdecomp}), we recognize that
\[\calP_{r}\Lambda^1(\square_2)\oplus\kappa\calH_r\Lambda^0(\square_2)\oplus d\calJ_{r+1}\Lambda^0(\square_2)=\calP_{r+1}^-\Lambda^1 \oplus d \calJ_{r+1} \Lambda^{0}(\square_2)=\calS_r^-\Lambda^1(\square_2),\]
since $\calJ_{r+1}\Lambda^1(\square_2)=\{0\}$.
For the third space, taking the $\flat$ operator for 2-forms, we have $\left(\mathscr{P}_r(x,y)\right)^\flat=\calP_{r}\Lambda^2(\square_2)=\calS_{r+1}^-\Lambda^2(\square_2)$.
This proves the first statement of Proposition~\ref{prop:cf-comp}.

The second statement of Proposition~\ref{prop:cf-comp} can be established similarly.
We state all the equivalencies first, then provide some details for the subtler cases:
\begin{align}
\left(\mathscr P_{r+1}\oplus\delta H^{3,I}_{r+1}\right)^\flat &= \calS_{r+1}^-\Lambda^0(\square_3) \\
\left(\pmb{\mathscr{P}}_{r}\oplus\bfx\times\tilde{\pmb{\mathscr{P}}}_{r}\oplus\nabla\delta H^{3,I}_{r+1}\oplus\delta E^{3,I}_{r+1}\right)^\flat &= \calS_{r+1}^-\Lambda^1(\square_3) \label{eq:cf-n3-k1}\\
\left(\pmb{\mathscr{P}}_{r}\oplus\bfx\tilde{\mathscr{P}}_{r}\oplus\nabla\times\delta E^{3,I}_{r+1}\right)^\flat &= \calS_{r+1}^-\Lambda^2(\square_3)  \label{eq:cf-n3-k2}\\
\left(\mathscr P_{r}\right)^\flat &= \calS_{r+1}^-\Lambda^3(\square_3)
\end{align}
The first and last statements are straightforward.
For the $1$-form case, (\ref{eq:cf-n3-k1}), we first recognize that
\[\left(\pmb{\mathscr{P}}_{r}\oplus\bfx\times\tilde{\pmb{\mathscr{P}}}_{r}\right)^\flat  = \calP_{r}^-\Lambda^1(\square_3). \]
We now claim that $\left(\nabla\delta H^{3,I}_{r+1}\right)^\flat =d\calJ_{r+1}\Lambda^0(\square_3)$.
It suffices to show that $\left(\delta H^{3,I}_{r+1}\right)^\flat =\calJ_{r+1}\Lambda^0(\square_3).$
The space $\delta H^{3,I}_{r+1}$ is defined as the span of polynomials of the form $xyz^{r+1}$ or $x\tilde{\mathscr{P}}_{r+1}(y,z)$, where  $\tilde{\mathscr{P}}_{r+1}(y,z)$ denotes homogeneous polynomial of degree $r+1$ in variables $y$ and $z$ only, or of similar forms with the variables permuted.
We have
\[\calJ_{r+1}\Lambda^0(\square_3) = \kappa\calH_{r+1,1}\Lambda^1(\square_3) \oplus \kappa\calH_{r+2,2}\Lambda^1(\square_3).\]
We can write $\calH_{r+1,1}\Lambda^1(\square_3)$ as the span of elements of the form $xp\,dy$ or $xp\,dz$ for any $p\in\tilde{\mathscr{P}}_{r}(y,z)$, or of similar forms with the variables permuted.
Observe that $\kappa\,xp\,dy = xyp$ and $\kappa\,xp\,dz = xzp$, both of which belong to $x\tilde{\mathscr{P}}_r(y,z)\subset\delta H^{3,I}_{r+1}$.
By similar analysis, after permuting variables, we have that $\left(\delta H^{3,I}_{r+1}\right)^\flat\supset \kappa\calH_{r+1,1}\Lambda^1(\square_3)$.
Next, the space $\calH_{r+2,2}\Lambda^1(\square_3)$ is spanned by the set $\{x^ryz\,dx$, $xy^rz\,dy, xyz^r\,dz\}$. 
Taking $\kappa$ of this set we get  $\{x^{r+1}yz, xy^{r+1}z, xyz^{r+1}\}$, establishing that $\left(\delta H^{3,I}_{r+1}\right)^\flat\supset \kappa\calH_{r+2,2}\Lambda^1(\square_3)$.
Since applying the flat operator to the elements of the spanning set for $\delta H^{3,I}_{r+1}$ produces a spanning set for $\calJ_{r+1}\Lambda^0(\square_3)$, we have established the claim.

Finally, we show that $\left(\delta E^{3,I}_{r+1}\right)^\flat=\calJ_{r+1}\Lambda^1(\square_3)$.
The space $\delta E^{3,I}_{r+1}$ is defined as the span of elements of the form $x\tilde{\mathscr{P}}_{r}(y,z)(y\nabla z - z\nabla y)$ or of two similar forms, with the variables permuted.
Let $p\in\tilde{\mathscr{P}}_{r}(y,z)$, i.e.\ $p$ is a homogeneous polynomial of degree $r$ in variables $y$ and $z$ only.
Observe that 
\[\left(xp (y\nabla z - z\nabla y)\right)^\flat= xp (y\,dz - z\,dy)=-\kappa (xp\,dydz)\in\kappa\calH_{r+1,1}\Lambda^2(\square_3).\]
Since $\kappa\calH_{r+1,1}\Lambda^2(\square_3)$ is spanned by form monomials that can be written as $\kappa (xp\,dydz)$ and similar form monomials with the variables permuted, we have that 
\[\left(\delta E^{3,I}_{r+1}\right)^\flat=\kappa\calH_{r+1,1}\Lambda^2(\square_3)=\calJ_{r+1}\Lambda^1(\square_3).\]
The last equality follows from (\ref{eq:jrk-def}), since any element of $\Lambda^2(\square_3)$ has linear degree at most 1.
By (\ref{eq:srm-directsumdecomp}), we have established (\ref{eq:cf-n3-k1}).
The final equality, (\ref{eq:cf-n3-k2}), can be confirmed by similar analysis.

\section{Proofs of dimension equality}
\label{appdx:dim-eq-exs}
We now prove that the number of degrees of freedom defined for the trimmed serendipity elements is equal to the dimension of the corresponding polynomial differential form space for $n=2$ and $n=3$.  
In our experience, all intuition for the cardinalities of these sets comes from the geometry of the $n$-cubes to which they are associated moreso than the algebra of binomial coefficients required for their computation.
Additional cases beyond those proved here can easily be checked using Mathematica or similar software, as we have done for $n=4$ and $n=5$ for $1\leq r\leq 100$.

Let $\dof(r,k,n)$ denote the number of degrees of freedom associated to $\calS_r^-\Lambda^k(\square_n)$; its value is defined by the formula (\ref{eq:Srm-dofs-count}).
Recall that $\dim\calS_r^-\Lambda^k(\square_n)$ can be computed using (\ref{eq:dim-Srmk}) and that $\dim\calJ_r\Lambda^k(\square_n)$ can be computed using Lemma~\ref{lem:dim-Jrk}.\\

\begin{remark}
In the following proofs, we adopt the convention:
\begin{equation}
\label{eq:binom-conven}
{n\choose k} :=\begin{cases}{\ds{n\choose k}} & \text{if}~n\geq k, \\ 0 & \text{if}~n<k. \end{cases}
\end{equation}
This convention is strictly for notational convenience as we frequently encounter summations whose upper index limit depends on $r$.  
For instance, the term $\ds{r-2\choose 2}$ appears in an expression for $\dof(r,0,2)$ only when $r\geq 4$.
By our convention, this summand is 0 when $r=1$, whereas converting it to the polynomial $\ds\frac{(r-2)(r-3)}{2}$ and evaluating at $r=1$ gives a value of 1.  
Hence, we will only convert binomial coefficients to functions when doing preserves the value according to the above convention. 
As we will see, this approach simplifies the presentation of the proofs.
\end{remark}

\newpage
\begin{prop}
For $k=0,1,2$, $\dof(r,k,2)=\dim\calS_r^-\Lambda^k(\square_2)$.
\end{prop}
\begin{proof}
We start with $k=0$. 
Expanding (\ref{eq:Srm-dofs-count}), we get
\[\dof(r,0,2)=4+4(r-1)+{r-2\choose 2}.\]
Note that our convention (\ref{eq:binom-conven}) applies to the third term in the sum, corresponding exactly to the summation index going from $d=0$ to $d=\min\{2,\lfloor r/2\rfloor\}$.
By (\ref{eq:srm-directsumdecomp}), we have that $\calS_r^-\Lambda^0(\square_2)=\calP_r^-\Lambda^0(\square_2)\oplus \calJ_r\Lambda^0(\square_2)$; the term $d\calJ_r\Lambda^{-1}$ is the empty set.
Using (\ref{eq:dim-prm}) and Lemma~\ref{lem:dim-Jrk}, we compute
\begin{align}
|\calP_r^-\Lambda^0(\square_2)| & = {r+2 \choose 2}, \notag\\
|\calJ_r\Lambda^0(\square_2)| & = 4+4(r-1)+{r-2\choose 2} - \sum_{j=0}^r(j+1) \notag \\
 & = 4+4(r-1)+{r-2\choose 2} - {r+2\choose 2}. \label{eq:dimJrL0S2}
\end{align}
Our convention (\ref{eq:binom-conven}) again applies to the term $\ds{r-2\choose 2}$ and again it appears in the summation only when $r\geq 4$.
Adding the formulae for $|\calJ_r\Lambda^0(\square_2)|$ and $|\calP_r^-\Lambda^0(\square_2)|$, we recover exactly the formula for $\dof(r,0,2)$.

Now we turn to $k=1$.
Expanding (\ref{eq:Srm-dofs-count}), we get
\[\dof(r,1,2)=\begin{cases} 4, & r=1, \\ r^2+2r+2, & r\geq 2. \end{cases}\]
The cases are due to the fact that $\min\{2,\lfloor r/2\rfloor+1\}=1$, for $r=1$, and $2$, for $r\geq 2$.
We have $\calS_r^-\Lambda^1(\square_2)=\calP_r^-\Lambda^1(\square_2)\oplus \calJ_r\Lambda^1(\square_2)\oplus d\calJ_r\Lambda^0(\square_2)$.
Recall that $|d\calJ_r\Lambda^0(\square_2)|=|\calJ_r\Lambda^0(\square_2)|$, for which we already have a formula from the $k=0$ case.
Again using (\ref{eq:dim-prm}) and Lemma~\ref{lem:dim-Jrk}, we can compute
\begin{align}
|\calP_r^-\Lambda^1(\square_2)| & = (r+2)r, \label{eq:dimPrmL1S2} \\
|\calJ_r\Lambda^1(\square_2)| & = {r+3\choose 2}-2{r+2\choose 2}+2{r\choose 2}-{r-1\choose 2}.\notag
\end{align}
Similar to the $k=0$ case, the term $\ds 2{r\choose 2}$ only appears for $r\geq 2$ and the term $\ds -{r-1\choose 2}$ only appears for $r\geq 3$, in accordance with our convention (\ref{eq:binom-conven}).
Converting the binomial coefficients to functions is still valid for any $r\geq 1$.  By doing so and simplifying, we find that
\begin{equation}
\label{eq:dimJrL1S2}
|\calJ_r\Lambda^1(\square_2)| =0.
\end{equation}
This was expected, given the general fact pointed out in~\cite[Equation (15)]{AA2014} that
\[\calJ_r\Lambda^k(\bbR^n)=0,\quad\text{for $k=n$ or $k=n-1$}.\]
Thus, we just have to add (\ref{eq:dimJrL0S2}) and (\ref{eq:dimPrmL1S2}) to compute the dimension of $\calS_r^-\Lambda^1(\square_2)$.
When $r=1$, we find that $\dim\calS_1^-\Lambda^1(\square_2)=4$.
For $r\geq 2$, converting to functions is valid; doing this and simplifying yields $r^2+2r+2$, recovering the formula for $\dof(r,1,2)$.

Finally we turn to the case $k=2$.
Expanding (\ref{eq:Srm-dofs-count}), we get
\[\dof(r,2,2)= {r+1\choose r-1} = {r+1 \choose 2}.\]
Note that $\calJ_r\Lambda^2(\square_2)=0$, since a 2-form in $n=2$ cannot be an image of $\kappa$, so $\calS_r^-\Lambda^2(\square_2) = \calP_r^-\Lambda^2(\square_2)\oplus\calJ_r\Lambda^1(\square_2)$.
We have $|\calP_r^-\Lambda^2(\square_2)| = \ds{r+1\choose 2}$ from (\ref{eq:dim-prm}) and $|\calJ_r\Lambda^1(\square_2)|=0$ from (\ref{eq:dimJrL1S2}), completing the proof.
\end{proof}
\begin{prop}
For $k=0,1,2,3$, $\dof(r,k,3)=\dim\calS_r^-\Lambda^k(\square_3)$.
\end{prop}
\begin{proof}
We begin by again recalling the observation from~\cite[Equation (15)]{AA2014} that
\[\calJ_r\Lambda^k(\bbR^n)=0,\quad\text{for $k=n$ or $k=n-1$}.\]
This allows the following simplifications of the decomposition (\ref{eq:srm-directsumdecomp}) for the spaces in $n=3$:
\begin{align*}
\calS_r^-\Lambda^0(\square_3) & = \calP_r^-\Lambda^0(\square_3)\oplus \calJ_r\Lambda^0(\square_3), \\
\calS_r^-\Lambda^1(\square_3) & = \calP_r^-\Lambda^1(\square_3) \oplus \calJ_r\Lambda^1(\square_3) \oplus d \calJ_r \Lambda^0(\square_3), \\
\calS_r^-\Lambda^2(\square_3) & = \calP_r^-\Lambda^2(\square_3)\oplus d \calJ_{r} \Lambda^{1}(\square_3),  \\
\calS_r^-\Lambda^3(\square_3) & = \calP_r^-\Lambda^3(\square_3).
\end{align*}
Since $|d \calJ_r \Lambda^k(\square_3)|=|\calJ_r \Lambda^k(\square_3)|$, we only need to expand the formula from Lemma~\ref{lem:dim-Jrk} for $\dim\calJ_r \Lambda^k(\square_3)$ for $k=0$ and $k=1$ to be able the compute the dimensions of all the spaces.
We find that
\begin{equation*}
|\calJ_r\Lambda^0(\square_3)| = \begin{cases} 4 & \text{if}~ r=1, \\ 10 & \text{if}~r=2, \\ 3(r+1) & \text{if}~r\geq 3. \end{cases}
\end{equation*}
and 
\begin{equation*}
|\calJ_r\Lambda^1(\square_3)| = \begin{cases} 2 & \text{if}~r=1, \\  3r & \text{if}~r\geq 2. \end{cases}
\end{equation*}
Using the formula for $|\calP_r^-\Lambda^k(\square_3)|$ from (\ref{eq:dim-prm}) and the above, we can produce formulae for $|\calS_r^-\Lambda^k(\square_3)|$ for each $k$.
We write out the formulae for $k=1$ in detail as it is the most elaborate:
\begin{equation}
\label{eq:dim-srm13}
|\calS_r^-\Lambda^1(\square_3)| = \begin{cases} 12,~36, & \text{if}~r=1,~2,~\text{respectively}, \\[2mm] \ds r{r+3\choose 2} + 3r + 3(r+1) & \text{if}~r\geq 3 \end{cases}
\end{equation}
Similarly, we can compute the degree of freedom count using (\ref{eq:Srm-dofs-count}) to produce formulae for $\dof(r,k,3)$ for each $k$.
In the case $k=1$, we get
\begin{equation}
\label{eq:dof-srm13}
\dof(r,1,3) = \begin{cases} 12,~36, & \text{if}~r=1,~2,~\text{resp.}, \\[2mm] \ds 6r^2 + 12 + (r - 2){r-3\choose 2} + (r - 1)(r - 3) & \text{if}~r\geq 3. 
\end{cases}
\end{equation}
In the case $r\geq 3$, we convert the binomial coefficients in (\ref{eq:dim-srm13}) and (\ref{eq:dof-srm13}) into polynomials and simplify, producing $\ds\frac{r^3}{2} + \frac{5 r^2}{2} + 9r + 3$ from each, thereby confirming the equality of the dimensions.
The remaining cases $k=0$, $k=2$, and $k=3$ are confirmed similarly.
\end{proof}

\end{document}